\documentclass[a4paper]{amsart} 
\usepackage[ascii]{inputenc} 

\usepackage{amssymb}
\usepackage{mathrsfs}
\usepackage{mathtools}

\newcommand{\Cal}{\mathcal}   
\newcommand{\itm}[1]{\item[{#1}]}   

\usepackage{xcolor}

\usepackage[shortlabels]{enumitem}
\setlist[enumerate,1]{label={(\Alph*)}}
\setlist[enumerate,2]{label={(\alph*)}}
\setlist[enumerate,3]{label={$\bullet_{\arabic*}$}}

\newenvironment{PROOF}[2][\proofname.]
   {\begin{proof}[#1]\renewcommand{\qedsymbol}{\textsquare$_{\mathrm{#2}}$}}
   {\end{proof}}

\newtheorem{theorem}{Theorem}[section]
\newtheorem{Theorem}[theorem]{Theorem}

\newtheorem{claim}[theorem]{Claim}
\newtheorem{Claim}[theorem]{Claim}
\newtheorem{conclusion}[theorem]{Conclusion}

\newtheorem{Corollary}[theorem]{Corollary}

\newtheorem{Lemma}[theorem]{Lemma}
\newtheorem{observation}[theorem]{Observation}
\newtheorem{Observation}[theorem]{Observation}

\theoremstyle{definition}

\newtheorem{definition}[theorem]{Definition}
\newtheorem{Definition}[theorem]{Definition}

\newtheorem{Example}[theorem]{Example}

\newtheorem{fact}[theorem]{Fact}
\newtheorem{Fact}[theorem]{Fact}
\newtheorem{hypothesis}[theorem]{Hypothesis}

\theoremstyle{remark}

\newtheorem{notation}[theorem]{Notation}

\newtheorem{remark}[theorem]{Remark}
\newtheorem{Remark}[theorem]{Remark}

\newtheorem{Criterion}[theorem]{Criterion}
\newtheorem{Assumption}[theorem]{Assumption}

\DeclareMathOperator{\Mod}{Mod}
\DeclareMathOperator{\Cat}{Cat}
\DeclareMathOperator{\Op}{Op}
\DeclareMathOperator{\EM}{EM}
\DeclareMathOperator{\GEM}{GEM}
\DeclareMathOperator{\sk}{sk}

\newcommand{\at}{\mathrm{at}}

\newcommand{\cf}{\mathrm{cf}}

\newcommand{\EC}{\mathrm{EC}}

\newcommand{\eq}{\mathrm{eq}}

\newcommand{\id}{\mathrm{id}}

\newcommand{\LST}{\mathrm{LST}}
\newcommand{\LS}{\mathrm{LS}}

\newcommand{\otp}{\mathrm{otp}}

\newcommand{\rang}{\mathrm{rang}}

\newcommand{\rng}{\mathrm{rng}}

\newcommand{\tp}{\mathrm{tp}}

\newcommand{\val}{\mathrm{val}}

\newcommand{\bfa}{\mathbf{a}}

\newcommand{\bfD}{\mathbf{D}}

\newcommand{\bfI}{\mathbf{I}}

\newcommand{\bfu}{\mathbf{u}}

\newcommand{\bbL}{\mathbb{L}}

\newcommand{\cP}{\mathscr{P}}

\newcommand{\gK}{\mathfrak{K}}

\newcommand{\gk}{\mathfrak{k}}

\newcommand{\varp}{\varepsilon}

\newcommand{\rest}{\restriction}

\newcount\skewfactor
\def\mathunderaccent#1#2 {\let\theaccent#1\skewfactor#2
\mathpalette\putaccentunder}
\def\putaccentunder#1#2{\oalign{$#1#2$\crcr\hidewidth
\vbox to.2ex{\hbox{$#1\skew\skewfactor\theaccent{}$}\vss}\hidewidth}}
\def\name{\mathunderaccent\tilde-3 }

\newbox\noforkbox \newdimen\forklinewidth
\setbox0\hbox{$\textstyle\bigcup$}
\setbox1\hbox to \wd0{\hfil\vrule width \forklinewidth depth \dp0
                        height \ht0 \hfil}
\setbox\noforkbox\hbox{\box1\box0\relax}
\def\unionstick{\mathop{\copy\noforkbox}\limits}
\def\nonfork#1#2_#3{#1\unionstick_{\textstyle #3}#2}
\def\nonforkin#1#2_#3^#4{#1\unionstick_{\textstyle #3}^{\textstyle
    #4}#2}
\setbox0\hbox{$\textstyle\bigcup$}
\setbox1\hbox to \wd0{\hfil{\sl /\/}\hfil}
\setbox2\hbox to \wd0{\hfil\vrule height \ht0 depth \dp0 width
                                \forklinewidth\hfil}
\newbox\doesforkbox
\setbox\doesforkbox\hbox{\box1\box0\relax}
\def\nunionstick{\mathop{\copy\doesforkbox}\limits}

\def\fork#1#2_#3{#1\nunionstick_{\textstyle #3}#2}
\def\forkin#1#2_#3^#4{#1\nunionstick_{\textstyle #3}^{\textstyle
    #4}#2}

\newcommand{\stickT}{%
\setbox255=\hbox{\raise1ex\hbox{$\hspace{0.2pt}\,\bullet\,$}}
\mathord{\rlap{\hbox to\wd255{\hss\hbox{$|$}\hss}}
\box255}
}
\newcommand{\stickS}{%
\setbox255=\hbox{\raise0.6ex\hbox{$\scriptstyle\bullet$}}
\mathord{\rlap{\hbox to\wd255{\hss\hbox{$\scriptstyle|$}\hss}}
\box255}
}

\usepackage[hidelinks]{hyperref}

\author[S. Shelah]{Saharon Shelah}

\address{Einstein Institute of Mathematics,
The Hebrew University of Jerusalem,
9190401, Jerusalem, Israel; and\\
Department of Mathematics,
Rutgers University,
Piscataway, NJ 08854-8019, USA}

\urladdr{https://shelah.logic.at/}

\thanks{First typed in the nineties.
The author would like to thank ISF-BSF for partially supporting  this research by a grant with Maryanthe  
Malliaris  
number NSF 2051825, BSF 3013005232. References like [Sh:950, Th0.2=Ly5] mean that the internal label of Th0.2 is y5 in Sh:950. The reader should note that the version on my website is usually more up-to-date than the one in arXiv.}

\makeatletter
\@namedef{subjclassname@2020}{\textup{2020} Mathematics Subject Classification}

\makeatother


\date{February 23, 2024} 

\newcommand\pmf{\par\medskip\par} 
\newcommand\un{\underbar}
\newcommand\e{\varepsilon}
\newcommand\vp{\varphi}
\newcommand\pbf{\par\bigskip}
\newcommand\w{\omega}
\renewcommand\k{\kappa}

\newcommand\uhr{\upharpoonright}
\def\ov{\bar}
\def\lng{\langle}
\def\om{\omega}
\def\heb{HEB} 
\renewcommand\a{\aleph}
\newcommand\lam{\lambda}

\def\prenice{\underset {\rm{nice}}  \preceq }
\def\al{\alpha}
\def\b{\beta}

\def\eqdf{\mathop=\limits^{\hbox{\rm def}}}
\def\del{\delta}

\def\z{\zeta}

\newcommand\subnice{\underset{\rm{nice}}  \subseteq}
\def\rng{\rangle}
\def\lng{\langle}
\def\d{\delta}
\def\al{\alpha}
\def\b{\beta}
\def\unl{\underline}

\def\lam{\lambda}
\def\prenice{\underset{\rm{nice}} \leq} 
\def\th{\theta}
\def\ov{\overline}

\def\lam{\lambda}
\def\lng{\langle}
\def\rng{\rangle}
\def\gam{\gamma}

\def\he{HEB H} 
\def\vav{HEB V} 

\begin{document}
\makeatletter\def\shfiuwefootnote{\gdef\@thefnmark{}\@footnotetext}\makeatother\shfiuwefootnote{Version 2024-03-03. See \url{https://shelah.logic.at/papers/E102/} for possible updates.}

\begin{abstract}
    In the original version of this paper, we assume a theory $T$ in the logic $\bbL_{\kappa, \aleph_{0}}$ is categorical in a cardinal $\lambda > \kappa$, and $\kappa$ is a measurable cardinal. There we prove that the class of models of $T$ of cardinality $<\lambda$ (but $\geq |T|+\kappa$) has the amalgamation property 
    under a natural order; this is a step toward understanding the character of such classes of models.

    In this revised version we replaced the class of models of $T$ by $\gk$, an AEC (abstract elementary class) which has $\LST$-number ${<} \, \kappa,$ or at least which behaves nicely for ultra-powers by $\bfD,$ some normal ultra-filter on $\kappa$ or just $\kappa$-complete non-principal ultra-filters on $\kappa$.

    Presently sub-section \S1A deals with $T \subseteq \bbL_{\kappa^{+}, \aleph_{0}}$ (and so does a large part of the introduction and little in the rest of \S1), but otherwise, all is done in the context of AEC.

    We leave the original introduction adding in the end after the three stars.
\end{abstract}

\title{Categoricity and amalgamation
for AEC,  \\ and $ \kappa $ 
measurable 
}

\author{Oren Kolman}

\thanks{On the old versions (publication number 362) the author expresses gratitude for the partial support of the Binational Science Foundation in this research and thanks Simcha Kojman for her unstinting typing work. 
}

\thanks{In the new version, the author  is  grateful for the generous funding of typing services donated by a person who wishes to remain anonymous and would like to thank the typist for the careful and beautiful typing.
This is publication number E102 in the second author
list.  
}


\maketitle

\newpage

\setcounter{section}{-1}

\newpage

\section*{Annotated Content}

\S0   \quad Introduction,  pg. \pageref{0}.

\S1  \quad  Preliminaries,  pg. \pageref{1}.

\begin{enumerate}
    \item[${{}}$] [In \S1A we review materials on fragments ${\mathcal F}$ of $\bbL_{\kappa, \aleph_0}$ (including the theory $T$) and basic model-theoretic properties (Tarski-Vaught property and L.S.), and we define amalgamation. In \S1B we move to AEC $\gk = (K, \leq_{\gk})$ which is our main framework now, and spell out  the connection. In \S1C, D we deal with indiscernibles and E.M. models, then we deal with limit ultra-powers which are \underbar{suitable} (for
    $\bbL_{\kappa, \aleph_{0}}$ and for our AECs) and in particular ultra-limits. Next, we introduce a notion basic for this paper $M {\underset{\text{nice}} \leq }N$ if there is a $\leq_{\gk}$-embedding of $N$ into suitable ultra-limit of $M$ extending the canonical embbeding.]
\end{enumerate}

\S2 \quad The amalgamation property for regular categoricity, pg. \pageref{2}.

\begin{enumerate}
    \item[${{}}$] [We   get amalgamation in $(K_\lambda, \leq_{\gk})$ when one of the extensions is nice, see Claim \ref{X2.1}. We prove that if $\gk$ is categorical in the regular $\lambda > \LST_{\gk} + \kappa$, then $(K_{< \lambda}, \leq_{\gk})$ has the amalgamation property. For this, we show that nice extension (in $K_{<\lambda}$) preserves being a non-amalgamation basis. We also start investigating (in  Theorem \ref{X2.5}) the connection between extending the linear order $I$ and the model $\EM(I)$: $I \underset{\text{nice}}  \subseteq J \Rightarrow \EM(I) \underset {\text{nice}}  \leq \EM(J)$; and give sufficient condition for $I
    \underset{\text{nice}}  \subseteq J$ (in Criterion \ref{X2.6}). From this, we get in $K_\lambda$ a model such that any sub-model of a suitable expansion is a $\underset{\text{nice}}  \leq$-sub-model (in Fact \ref{X2.7}, Theorem \ref{X2.10}(2)), and conclude the amalgamation property in $(K_{<\lambda}, \leq_{\gk})$ when $\lambda$ is regular (in Theorem \ref{X2.9}) and something for singulars in Theorem \ref{X2.10}.]
\end{enumerate}

\S3 \quad Toward removing the assumption of regularity from the existence of universal extensions,  pg. \pageref{3}.  

\begin{enumerate}
    \item[${{}}$] [The problem is that $\EM(\lambda)$ has many sub-models which ``sit'' well in it and we can prove that there are many amalgamation bases but we need to get this simultaneously. First in Theorem \ref{X3.1} we show that, if $\langle M_i \colon i< \theta^+\rangle$ is $\leq_{\gk}$-increasing continuous sequence of models from $K_\theta,$ then for a club of $i< \theta^+$ we have $M_i \underset{\text{nice}}{\leq}  \bigcup \{M_j \colon  j< \theta^+\}$. In Definition \ref{X3.5}, we define nice models (essentially, every reasonable extension is nice). Next (in Theorem \ref{X3.3}) we show that nice models are dense in $K_{\theta}$. Also (by Theorem \ref{X3.4}) many embeddings are nice and (in Corollary \ref{X3.5}) we show that being nice implies being amalgamation base. Then we define a universal extension of $M\in K_\theta$ in $K_\partial$ (Definition \ref{X3.6}), we prove existence over a model in Lemma \ref{X3.7} and after preparation  prove the existence (Corollary \ref{X3.10}, Corollary \ref{X3.11}).]
\end{enumerate}

\S4 \quad $(\theta, \partial)$-saturated models,  pg. \pageref{4}.

\begin{enumerate}
    \item[${{}}$] [If $M_i \in K_\theta$ for $i\leq \partial$ is increasing continuous, $M_{i+1}$ universal over $M_i$, and each $M_i$ is nice, then\footnote{In \cite{Sh:600} we say $M_{j}$ is $(\theta, \partial)$-brimmed.} we say $M_\partial$ is
    \emph{$(\theta, \partial)$-saturated} over $M_0$. We show existence (and uniqueness). We connect this to more usual saturation and prove that $(\theta, \partial)$-saturation implies niceness (in Theorem \ref{X4.10}).]
\end{enumerate}

\S5 \quad  The amalgamation property for $K_{<\lambda}$, pg. \pageref{5}.

\begin{enumerate}
    \item[${{}}$] [After preliminaries we prove that for $\theta\leq \lambda$ (and $\theta \geq \LST(\gk) + \kappa$ of course) every member of $K_\theta$ can be extended to one with many nice sub-models, this is done by induction on $\theta$ using the niceness of $(\theta_1, \partial_1)$-saturated models. Lastly, we conclude that every $M\in K_{<\lambda}$ is nice hence $K_{<\lambda}$ has the amalgamation property.] 
\end{enumerate}

\newpage

\section{Introduction}\label{0}

The main result\footnote{In the old version.} of this paper is a proof of the following theorem:

\begin{Theorem}
    Suppose that $T$ is a theory in a fragment of $\bbL_{\kappa, \aleph_{0}}$ where $\k$ is a measurable cardinal.  If $T$ is categorical in the cardinal $\lambda > \k + | T | $, then  ${\mathcal K}_{< \lambda}$, the class of models of $T$ of power strictly less than $\lambda$ (but $\geq \chi = \kappa+ |T|$), has the amalgamation property (see Definition \ref{X1.3} (1)(2)).
\end{Theorem}

The interest in this theorem stems in part from its connection with the study of categoricity spectra.  For a theory $T$ in a logic $\mathscr{L}$ let us define $\Cat(T)$, the \emph{categoricity spectrum of $T$}, to be the collection of those cardinals $\lambda$ in which $T$ is categorical.  In the 1950's \L os
conjectured that if $T$ is a countable theory in first-order logic, then $\Cat (T)$ contains every uncountable cardinal or no uncountable cardinal.  This conjecture, based on the example of algebraically closed fields of fixed characteristic, was verified by Morley \cite{Mo65}, who proved that if a countable first-order theory is categorical in some uncountable cardinal, then it is categorical in every uncountable cardinal.  Following advances made by Rowbottom \cite{Row64}, Ressayre \cite{Res69} and Shelah \cite{Sh:2}, Shelah \cite{Sh:31} proved the \L os conjecture for uncountable first-order theories: if $T$ is a first-order theory categorical in some cardinal $\lambda > |T| + \a_0$, then $T$ is categorical in every cardinal $\lambda > |T| + \a_0$.  It is natural to ask whether analogous results hold for theories in logics other than first-order logic.  Perhaps the best-known extensions of first-order logic are the infinitary logics $\bbL_{\lambda, \kappa}$.  As regards theories in $\bbL_{\kappa, \aleph_{0}}$, Shelah (see \cite{Sh:87a} and \cite{Sh:87b})  continuing work begun in \cite{Sh:48} introduced the concept of excellent classes: these have models in all cardinalities, have the amalgamation property and satisfy the \L os conjecture.  In particular, if $\varphi$ is an excellent sentence of $\bbL_{\aleph_{1}, \aleph_{0}}$, then the \L os conjecture holds for $\varphi$.  Furthermore, under some set-theoretic assumptions (weaker than the Generalized Continuum Hypothesis) if $\varphi$ is a sentence in $\bbL_{\aleph_{1}, \aleph_{0}}$ which is categorical in $\a_n$ for every natural number $n$ (or even just if $\varphi$ is a sentence in $\bbL_{\kappa, \aleph_{0}}$ with at least one uncountable model not having too many models in each $\a_n$), then $\varphi$ is excellent. Now, \cite{Sh:300}, \cite{Sh:h} try to develop classification theory in some non-elementary classes. We cannot expect much for $\bbL_{ \lambda, \kappa}$ for $\kappa > \aleph_0$. The first author conjectured that if $\varphi$ is a sentence in $\bbL_{\aleph_{1}, \aleph_{0}}$ categorical in some $\lambda > \beth_{\omega_1}$, then $\varphi$ is categorical in every $\lambda > \beth_{\omega_1}$.  (Recall that the Hanf number of $\bbL_{\aleph_{1}, \aleph_{0}}$ is $\beth_{\omega_1}$, so if $\psi$ is a sentence in $\bbL_{\aleph_{1}, \aleph_{0}}$ and $\psi$ has a model of power $\lambda \ge \beth_{\omega_1}$, then $\psi$ has a model in every power $\lambda \ge \beth_{\omega_1},$ see  \cite{Ke71}) .  There were some who asked why so tardy the beginning.  Recent work of Hart and Shelah \cite{Sh:323} showed that for every natural number $k$ greater than 1 there is a sentence $\psi_k$ in $\bbL_{\aleph_{1}, \aleph_{0}}$ which is categorical in the cardinals $\a_0, \dots, \a_{k-1}$, but which has many models of power $\lambda$ for every cardinal $\lambda \ge 2^{\a_{k-1}}$.  The general conjecture for
$\bbL_{\aleph_{1}, \aleph_{0}}$ remains open nevertheless.  As regards theories in $\bbL_{\kappa, \aleph_{0}}$, progress has been recorded under the assumption that $\k$ is
a strongly compact cardinal.  Under this assumption Shelah and Makkai \cite{Sh:285}  have established the following results for a $\lambda$-categorical theory $T$ in a fragment ${\mathcal F}$ of $\bbL_{\kappa, \aleph_{0}}$:

1) if $\lambda$ is a successor
cardinal and $\lambda > (( \k')^\k)^+$ where $\k' = \max ( \k, | {\mathcal F}| )$, then $T$ is categorical in every cardinal greater than or equal to $\min (\lambda, \beth_{(2^{\k'})^+} )$,

2) if $\lambda > \beth_{\k+1} (\k')$, then
$T$ is categorical in every cardinal of the form $\beth_\delta$ with $\delta$ divisible by $(2^{\k'})^+$ (i.e. for some ordinal $\alpha > 0,$ $\delta = (2^{\kappa'})^+\cdot \alpha$ (ordinal multiplication)).

In proving theorems of this kind, one has recourse to the amalgamation property which makes possible the construction of analogs of saturated models.  In turn, these are
of major importance in categoricity arguments.  The amalgamation property holds for theories in first-order logic \cite{CK73} and in $\bbL_{\k, \kappa}$ when $\k$ is a strongly compact cardinal (see e.g. \cite{Sh:285}: although $\prec_{\bbL_{\k,\k}}$
fails the Tarski-Vaught property for unions of chains of length $\k$ (whereas $\prec_{\bbL_{\kappa, \aleph_{0}}}$ satisfies it), under a categoricity assumption it can be shown that $\prec_{\bbL_{\kappa, \aleph_{0}}}$ and $\prec_{\bbL_{\k, \k}}$ coincide). However, it is not known in general for theories in $\bbL_{\kappa, \aleph_{0}}$ or $\bbL_{\k, \k}$ when one weakens the assumption on $\k$, in particular when $\k$ is just a measurable cardinal.  Nevertheless, categoricity does imply the existence of reasonably saturated models in an appropriate sense, and it is possible to begin classification theory.  This is why the main theorem of the present paper is of relevance regarding the categoricity spectra of theories in $\bbL_{\kappa, \aleph_{0}}$ when $\k$ is measurable.
    
A sequel to this paper under preparation (which is now \cite{Sh:472}) tries to provide a characterization of $\Cat(T)$ at least parallel to that in \cite{Sh:285} and we hope to deal with the corresponding classification theory later.  This division of labor both respects historical precedent and is suggested by the increasing complexity of the material. Another sequel deals with abstract elementary classes (in the sense of \cite{Sh:88}) (see \cite{Sh:472}, \cite{Sh:394} respectively). On more work see \cite{Sh:576}, \cite{Sh:600}.
     
The paper is divided into five sections. Section 1 is preliminary and notational.  In section 2 it is shown that if the theory $T \subseteq \bbL_{\kappa, \aleph_{0}}$ or just suitable AEC $\gK$ is categorical in the regular cardinal $\lambda > \k + |T|$, then $K_{< \lambda}$ has the amalgamation property.  Section 3 deals with weakly universal models, section
4 with $(\theta, \partial)$-saturated and $\bar{\theta}$-saturated models.  In section 5 the amalgamation property for $K_{< \lambda}$ is established.
    
All the results in this paper (other than those explicitly credited) are due to Saharon Shelah. $$ \ast \ \  \ \ \ast \ \ \  \ \ast  $$

On a more recent survey see \cite{Sh:E56} and a recent one see \cite{Sh:842}, in particular on the history of $\kappa$-compact AEC. 

We had stated that clearly, the proof of \cite{Sh:362} works for AEC, but the referee of \cite{Sh:842} 
asked to do it explicitly. Here we justify \cite[4.7]{Sh:842}. Note that, \cite[1.1, ~1.2]{Sh:362} essentially proves that ($\Mod(T),\prec_{T}$) is an AEC ignoring Ax. V of AEC (see Definition \ref{u2}), so Fact \ref{X1.2}(2) was added. 

We thank Shimoni Garti for his help in proofreading and the referee for pointing out some obscure points. 
    
\newpage

\section{PRELIMINARIES}\label{1}

To start things off in this section, let us fix notation, provide basic definitions and well-known facts, and formulate our working assumptions.

The working assumptions in force throughout the paper are these.

\begin{Assumption}\label{X1.-3}
    The cardinal $\k$ is an uncountable measurable cardinal, and so there is a $\k$-complete non-principal ultra-filter on $\k,$ we can fix such $\bfD$.
\end{Assumption}

\begin{Assumption}\label{X1.-2}\ 

    (1) The theory $T$ is a theory in the infinitary logic $\bbL_{\kappa, \aleph_{0}}$, $\chi = \kappa + \vert T \vert$ and vocabulary $\tau = \tau_{T}$, \underline{or}

    (2) $\gk$ is an AEC which is $\bfD$-compact (see Definition \ref{u2} and Definition \ref{y1} respectively) and $\chi = \kappa + \LST(\gk).$
\end{Assumption}

Our main theorem for the logic $\bbL_{\kappa, \aleph_{0}}$ is: 

\begin{theorem}\label{a5}
    If $T \subseteq \bbL_{\kappa, \aleph_{0}}$ is categorical in $\lambda > \kappa + \vert \tau_{T} \vert$ \underline{then} the class of models of $T$ of cardinality $< \lambda$ but $\geq \kappa + \vert \tau_{T} \vert$ (under the so called $\prec_{\mathcal{F}},$ see Notation \ref{a8}(2)(7), (8)) has the amalgamation property.
\end{theorem}

\begin{PROOF}{\ref{a5}}
    Use Theorem \ref{u8} on AEC which is applicable by Conclusion \ref{u20} and recalling Definition \ref{u14} and Claim \ref{u17}.  
\end{PROOF}

From these assumptions follow certain facts, of which the most important are these.

\begin{fact}\label{X1.-1}
    For each model $M$ of $T$, $\k$-complete ultra-filter $D$ over $I$ and suitable set $G$ of equivalence relations on $I \times I$ (see Definition \ref{w11}) the limit ultra-power $\Op(M) = \Op(M,I,D,G)$ is a model of $T$.
\end{fact}

\begin{fact}\label{X1.0}
    For each linear order $I = (I, \le)$ there exists an Ehrenfeucht-Mostowski model $\EM(I)$ of $T$ (see Definition \ref{v0}(6)). 
\end{fact} 

This section is divided into several subsections: in \S1A we deal with a theory $T$ in $\bbL_{\kappa, \aleph_{0}},$ in \S1B we move to AEC $\gk$ showing that the context in \S1A is a special case. Then in \S1C we deal with $\EM$ models. Finally, in \S1D we deal with ultra-powers, ultra-limits, and nice sub-models.  

\subsection{Frame for $\bbL_{\kappa, \aleph_{0}}$}\label{1A}\ 

Relevant set-theoretic and model-theoretic information on measurable cardinals can be found in \cite{Je03}, \cite{CK73}, and \cite{Dic75}. 

\begin{notation}\label{a1}
    Let $\tau$ denote\footnote{In the old version it was called ``language'' and denoted by $L$.} a vocabulary, i.e. a set of finitary relation and function symbols, including equality (i.e. the arity of the symbols in $\tau_{\sk}$ is always finite). So  $|\tau|$ is the cardinality of the vocabulary $\tau$.
\end{notation}

\begin{definition}\label{a2}\ 

    (1) For cardinals $\kappa \le \lambda, \  \bbL_{ \lambda, \kappa}$ is the logic such that for any vocabulary $\tau,$ $\bbL_{\lambda, \kappa}(\tau)$ is the smallest set of (possibly infinitary) formulas in the vocabulary $\tau$ which contains all first-order formulas and which is closed under: 

    \begin{enumerate}
        \item[(A)] the formation of conjunctions (disjunctions) of any set of formulas of power less than $\lambda$, provided that the set of free variables in the conjunctions (disjunctions) has power less than $\kappa$,

        \item[(B)] the formation of $\forall \bar x \vp,  \exists \bar x \vp$, where $\bar x = \langle x_\alpha \colon \alpha < \alpha_{\ast} \rangle$ is a sequence of variables of length $\alpha_{\ast} < \kappa$. 
    \end{enumerate}

    (2) Whenever we use the notation $\varphi(\bar{x})$ to denote a formula in $\bbL_{\lambda, \kappa},$ we mean that $\bar{x}$ is a sequence $\langle x_{\alpha} \colon \alpha < \alpha_{\ast} \rangle$ as above. So if $\varphi(\bar{x})$ is a formula in $\bbL_{\kappa, \aleph_{0}},$ then $\bar{x}$ is a finite sequence of variables. 

    (3) So $\bbL = \bbL_{\aleph_{0}, \aleph_{0}}$ is a first order logic. 
\end{definition}

\begin{notation}\label{a8}\ 

    (1) $\Cal F$ denotes a fragment of $\bbL_{\kappa, \aleph_{0}}(\tau)$, i.e. a set of formulas of $\bbL_{\kappa, \aleph_{0}}(\tau)$ which contains all atomic formulas of $\tau$, and which is closed under negations, finite conjunctions (finite disjunctions), and the formation of subformulas.  An $\Cal F$-formula is just an element of $\Cal F$.
    
    (2) $T$ is a theory in $\bbL_{\kappa, \aleph_{0}}(\tau)$, so there is a fragment $\Cal F$ of $\bbL_{\kappa, \aleph_{0}}$ such that $T \subseteq \Cal F$ and $| \Cal F | < | T|^+ + \k$. Let $\mathcal{F}_{T}$ be the minimal such $\mathcal{F}$. If not said otherwise, $T$ and $\mathcal{F} = \mathcal{F}_{T}$ are fixed. 
        
    (3) Models of $T$ (invariably referred to as models) are $\tau$-structures which satisfy the sentences of $T$. They are generally denoted $M, N, \dots, $ and $  |M|$ is the universe of the $\tau$-structure $M; \; \Vert M \Vert$ is the cardinality of $|M|$.
    
    (4) For a set $A, \; |A|$ is the cardinality of $A$ and $^{< \omega} A$ is the set of finite sequences in $A$ and for $\bar a = \langle a_0 \dots a_{n-1} \rangle \in \; ^{< \omega} A, \ \rm{lg} (\bar a) = n$ is the length of $\bar a$. Similarly, if $\bar a=\lng a_\z:\z<\del\rng$, we write $\lg(\bar a)=\del$, where $\del$ is an ordinal.
    
    (5) For an element $R$ of $\tau$ and a $\tau$-model $M$, let $\val(M, R)$, or $R^M$, be the interpretation of $R$ in the $\tau$-structure $M$.
        
    (6) We ignore models of power less than $\k$.  $K$ is the class of all models of $T$; $$ K_\lambda = \{ M \in K: || M || = \lambda \}, \;  K_{< \lambda} = \underset{\mu < \lambda} \bigcup K_\mu, \; K_{\le \lambda} = \underset{\mu \le \lambda} \bigcup K_\mu, \, K_{[\mu, \lambda)} = \underset{\mu \le \chi< \lambda} \bigcup K_\chi. $$
        
    (7) We write $f  \colon M \underset{\Cal F} \rightarrow N$ (may be abbreviated $f: M \rightarrow N$) to mean
    that $f$ is an $\Cal F$-elementary embedding (briefly, an embedding) of $M$ into $N$,  i.e. $f$ is a function with domain $|M|$ into $|N|$ such that for every $\Cal F$-formula $\vp (\bar x)$, and $\bar a \in \; ^{< \omega} |M|$ with $\lg(\bar a) = \lg (\bar x),  M \vDash \vp [\bar a]$ iff $N \vDash \vp [f (\bar a)]$, where if $\bar a = \langle a_i: i < n \rangle$, then $f (\bar a) := \langle f (a_i): i < n \rangle$.
    
    (8) In the special case where an embedding $f$ is a set-inclusion (so that $|M| \subseteq  \vert N \vert$), we write $M \prec_{\Cal F} N$ (briefly $M \prec N$), instead of $f: M \underset{\Cal F} \rightarrow N.$ We may say that $M$ is an $\Cal F$-elementary sub-model of $N$, or $N$ is an $\Cal F$-elementary extension of $M$.
\end{notation}
    
\begin{notation}\label{a11}\ 

    (1) $(I, \le_I), (J, \le_J)$ are partial orders; we will not bother to subscript the order relation unless really necessary; we may write $I$ for $(I, \le).$ We say $(I, \le)$ is directed iff for every $i_1$ and $i_2$ in $I$, there is $i \in I$ such that $i_1 \le i$ and $i_2 \le i$. $(I, <)^\ast$ is the (reverse) partial order $(I^\ast, <^\ast)$ where $I^\ast = I$ and $s <^\ast t$ iff $t < s$.
        
    (2) A sequence $\langle M_i: i \in I \rangle$ of models indexed by $I$ is a $\prec_{\Cal F}$-directed system iff $(I, \le)$ is a directed partial order and for $i \le j$  in $I, M_i \prec_{\Cal F} M_j$.

    Note that, the union $\underset{i \in I} \cup M_i$ of a $\prec_{\Cal F}$-directed system $\langle M_i  \colon i \in I \rangle$ of $\tau$-structures is an $\tau$-structure. In fact, more is true.
\end{notation}

\begin{Fact}\label{X1.1}\  

    (1) (Tarski-Vaught property) The union of a $\prec_{\Cal F}$-directed system $\langle M_i : i \in I \rangle$ of models of $T$ is a model of $T$, and for every $j \in I, M_j \prec_{\Cal F} \underset{i \in I} \cup M_i$.
    
    (2) For $\bar{M}$ as above, if $M$ is a fixed model of $T$ such that for every $i \in I$ there is $f_i \colon M_i \underset{\Cal F} \rightarrow M,$ $I$ is directed, and for all $i \le j$ in $I, f_i \subseteq f_j$, then $\underset{i \in I} \cup f_i: \underset{i \in I} \cup M_i \underset{\Cal F} \rightarrow M$. In particular, if $M_i \prec_{\Cal F} M$ and $f_{i}$ is the identity function on $M_{i}$ for every $i \in I$, then $\underset{i \in I}  \cup M_i \prec_{\Cal F} M$. Let $\alpha$ be an ordinal.  A $\prec_{\Cal F}$-chain of models of length $\alpha$ is a sequence $\langle M_\beta: \beta < \alpha \rangle$ of models such that if $\beta < \gamma < \alpha$, then $M_\beta \prec_{\Cal F} M_\gamma$.  The chain is continuous if for every limit ordinal $\beta < \alpha,  M_\beta = \underset{\gamma < \beta}  \cup M_\gamma$.
\end{Fact}

\begin{Fact}\label{X1.2}\      
    
    (1) (Downward L\"owenheim-Skolem Property):
    Suppose that $M$ is a model of $T$, $ A \subseteq | M |$ and $\max (\k + | T |, |A| ) \le \lambda \le || M ||$.  \underline{Then} there is a model $N$ such that $A \subseteq | N |,  || N || = \lambda$ and $N \prec_{\mathcal{F}} M$.
    
    (2) If $N$ and $M_{1} \subseteq M_{2}$ are $\tau$-models, $\mathcal{F}$ is a fragment of $\bbL_{\kappa, \aleph_{0}},$ and $M_{\ell} \prec_{\mathcal{F}} N$ for $\ell = 1, 2$ \underline{then} $M_{1} \prec_{\mathcal{F}} M_{2}.$
\end{Fact}



Now we turn from the rather standard model-theoretic background to the more specific concepts which are central in our investigation.

\begin{Definition}\label{X1.3}\ 

    (1) Suppose that $<$ is a binary relation on a class $K$ of models (mainly $(K, <) = (K_{\gk}, <_{\gk}$), see below). We say $\Cal K = \langle K, < \rangle$ has the \emph{amalgamation property} (AP) \underline{iff} for every $M, M_1, M_2 \in K$, if $f_i$ is an isomorphism from $M$ onto $\mathrm{rng}(f_i)$ and $\mathrm{rng}(f_i)<M_i$ for $i=1,2$, then there exist $N\in K$ and isomorphisms $g_i$ from $M_i$ onto $\mathrm{rng}(g_i)$ for $i = 1,2$ such that $\mathrm{rng}(g_i) < N$ and $g_1 f_1 = g_2 f_2$. The model $N$ is called an \emph{amalgam} of $M_1, M_2$ over $M$ with respect to $f_1, f_2$.
    
    (2) An $\tau$-structure $M$ is an \emph{amalgamation base} (a.b.) for $\Cal K = \langle K, < \rangle$ iff $M \in K$ and whenever for $i = 1, 2, M_i \in K$ and $f_i$ is an isomorphism from $M$ onto $\mathrm{rng}(f_i), \, \mathrm{rng}(f_i) < M_i$, then there exist $N \in K$ and isomorphisms $g_i \ (i = 1,2)$ from $M_i$ onto $\mathrm{rng}(g_i)$ such that $\mathrm{rng}(g_i) < N$ and $g_1 f_1 = g_2 f_2$.

    (3) We say  $\Cal K = \langle K, < \rangle$ has AP iff every model in $K$ is an a.b. for $\Cal K$.
\end{Definition}

\begin{Example}\label{X1.3A}
    Suppose that $T$ is a theory in first-order logic having an infinite  model.  Define, for $M, N$ in the class $K_{\le | T| + \a_0}$ of models of $T$ of power at most $|T|+\a_0$, $M < N$ iff  the indentity on $\vert M \vert$ is an embedding of $M$ onto an elementary sub-model of $N.$   Then $\Cal K_{\le |T| + \a_0} = \langle K_{\le |T| + \a_0}, < \rangle$ has AP, (see \cite{CK73}).
\end{Example}

\begin{Example}\label{X1.3B}
    Suppose that $T$ is a theory in $\bbL_{\kappa, \aleph_{0}}$ and $\Cal F$ is  a fragment of $\bbL_{\k, \aleph_{0}}$ containing $T$ with $| \Cal F| < |T|^+ + \k$. Let $<$ be the binary relation $\prec_{\Cal F}$ defined on the class $K$ of all models of $T$.  $M \in K$ is an a.b. for $\Cal K$ iff whenever for $i = 1, 2, M_i \in K$ and $f_i$ is an $\prec_{\Cal F}$-elementary embedding of $M$ into $M_i$, there exist $N \in K$ and $\Cal F$-elementary embeddings $g_i (i = 1,2)$ of $M_i$ into $N$ such that $g_1 f_1 = g_2 f_2$.
\end{Example}

\begin{Definition}\label{X1.4}
    Suppose that $<$ is a binary relation on a class $K$ of models. Let $\mu$ be a cardinal.
    $M \in K_{\le \mu}$ is a \emph{$\mu$-counter amalgamation basis} $(\mu$-c.a.b.) of
    $\Cal K = \langle K, < \rangle$ iff
    there are $M_1, M_2 \in K_{\le \mu}$ and isomorphisms $f_i$ from $M$ into
    $M_i$ such that:
    
    \begin{enumerate}
        \item[(a)] $\mathrm{rng}(f_i) < M_i (i = 1,2)$,
    
        \item[(b)] there is no amalgam $N \in K_{\le \mu}$ of $M_1, M_2$ over $M$ with respect to $f_1, f_2$.
    \end{enumerate}
\end{Definition}

\begin{Observation}\label{X1.5}
    Suppose that $T, \Cal F$ and $<$ are as in Example \ref{X1.3B} and $\k + |T| \le \mu < \lambda$.  Note that if there is an amalgam $N'$ of $M_1, M_2$ over $M$ (for $M_1, M_2, M$ in $K_{\le \mu}$), then  by Fact \ref{X1.2} there is an amalgam $N \in K_{\le \mu}$ of $M_1, M_2$ over $M$.
\end{Observation}


\subsection{Replacing $T$ by AEC}\label{B}\ 

On AEC see \cite{Sh:88}, \cite{Sh:88r} or \cite{Bal09}, recall:

\begin{definition}\label{u2}
    We say $\gk = (K_{\gk}, \leq_{\gk})$ is an a.e.c. with 
    L.S.T. number $\lambda(\gk) = \LS_{\gk} = \LST(\gk)$, we may write $K$ for $K_{\gk},$ \underline{when} $K$ is a class of $\tau_{\gk}$-models, $\leq_{\gk}$ a two-place relation on $K$ and 

    \begin{itemize}
        \item     \underline{Ax 0}:  The holding of $M \in K,N \le_\gk M$  depend on $N,M$ only up to isomorphism, i.e. $[M \in K,M \cong N \Rightarrow N \in K]$ and 
        [if $N \le_\gk M$ and $f$ is an isomorphism from $M$ onto the $\tau$-model $M'$ and $ f \, \rest \,  N$ is an isomorphism from $N$ onto $N'$ \underline{then} $N' \le_\gk M'$.]

        \item     \underline{Ax I}:  if $M \le_{\gk} N$  then $M \subseteq N$ (i.e. $M$ is a sub-model of $N$).

        \item     \underline{Ax II}:  $M_0 \le_{\gk} M_1 \le_{\gk} M_2$ implies $M_0 \le_{\gk} M_2$ and $M \le_{\gk} M$ for $M \in K$.

        \item     \underline{Ax III}:  If $\lambda$ is a regular cardinal, $M_i \ (i < \lambda)$ is a $\le_{\gk}$-increasing (i.e. $i < j < \lambda$ implies  $M_i \le_\gk M_j$) and continuous (i.e. for every limit ordinal $\delta < \lambda,M_\delta = \bigcup_{i < \delta} M_i$) \underline{then} $M_0 \le_{\gk} \bigcup_{i < \lambda} M_i$. Hence $M_{j} \leq_{\k} \bigcup_{i < \lambda} M_{i}$ for every $j < \lambda.$  

        \item     \underline{Ax IV}:  If $\lambda$ is a regular cardinal and $M_i$  (for $i < \lambda)$ is $\le_\gk$-increasing continuous and $M_i \le_\gk N$ for $i < \lambda$ 
        \underline{then}  $\bigcup_{i < \lambda} M_i \le_{\mathfrak{k}} N$.

        \item     \underline{Ax V}:  If $N_0 \subseteq N_1 \le_\gk M$ and $N_0 \le_{\mathfrak{k}} M$  \underline{then} $N_0 \le_\gk N_1$.

        \item     \underline{Ax VI}:  If $A \subseteq N \in K$ and $|A| \le \text{ LST}(\gk)$  then for some $M \le_\gk N,A \subseteq |M|$ and $\|M\| \le  \text{ LST}(\gk)$ (and $\text{\rm LST}(\gk)$ is the minimal infinite cardinal satisfying this axiom which is $\ge |\tau|$;  the $\ge |\tau|$ is for notational simplicity). 
    \end{itemize}
\end{definition}

\begin{definition}\label{u5}\ 

    (1) We define ``$\gk$ categorical in $\lambda$'', $\gk_{< \lambda}$, ``$\gk$ has amalgamation'' ``$M \in K_{\gk}$ is a.b.'', ``$M$ is c.a.b.'' naturally (see Definitions \ref{X1.3} and \ref{X1.4}). 

    (2) Let $\gk_{\lambda} = (K_{\lambda}, \leq_{\gk} \rest K_{\lambda}),$ where $K_{\lambda} = \{ M \in K_{\gk} \colon \Vert M \Vert = \lambda \}.$

    (3) For $\chi < \lambda,$ let $\gk_{[\chi, \lambda)} = (K_{[\chi, \lambda)}, \leq_{\gk} \rest K_{[\chi, < \lambda)]}),$ where $K_{[\chi, \lambda)} = \bigcup \{ K_{\mu} \colon \mu \in [\chi, \lambda) \}.$ 
\end{definition}

So our main theorem is: 

\begin{Theorem}\label{u8}
    Assume $\kappa$ is a measurable cardinal, $\gk$ is an AEC, and $\chi = \LST_{\gk} + \kappa < \lambda$, and $\LS_{\gk} < \kappa$ or just $\gk$ is $\bfD$-compact (see Definition \ref{X1.-2} and Assumption \ref{X1.-3}).  If $\gk$ is categorical in $\lambda$ \underline{then} $\gk_{[\chi, \lambda)}$ has amalgamation, see Definition \ref{X1.3}. 
\end{Theorem}

\begin{PROOF}{\ref{u8}}
    First, without loss of generality, assume that Hypothesis \ref{b2} holds.

    [Why? If $\LS_{\gk} < \kappa$ then by Claim \ref{v8}(0), without loss of generality $\vert \tau_{\gk} \vert \leq 2^{\LS(\gk)}$, hence $\vert \tau_{\gk} \vert < \chi$ and by Claim \ref{y3}(1), $\gk$ is $\bfD$-compact (see Assumption \ref{X1.-3}). So in any case $\gk$ is $\bfD$-compact and by Claim \ref{b0}, Hypothesis \ref{b2}(1) holds.

    By Claim \ref{v8}(1), (2) also Hypothesis \ref{b2}(2) holds. So Hypothesis \ref{b2} holds indeed.] 

    Recall that, in \S2-\S5 we assume Hypothesis \ref{b2}. 

    Second, if $\lambda$ is regular, then the desired conclusion holds by \S2, that is, by Theorem \ref{X2.9}. 

    Third, if $\lambda$ is singular, then the desired conclusion holds by \S5, that is, by Corollary \ref{X5.5}. 
\end{PROOF}

\begin{claim}\label{u11}
    Assume $\gk$ is an AEC and $\tau = \tau_{\gk}.$ Then;

    There are $\tau_{1} = \tau_{\gk, 1} \supseteq \tau_{\gk}$ of cardinality $\vert \tau \vert + \LST_{\gk}$ and a set $\cP$ of q.f. (quantifier free) $1$-types in $\bbL(\tau_{1})$ such that: 

    \begin{enumerate}
        \item[(A)] a $\tau$-structure $M$ belongs to $K_{\gk}$ \underline{iff} it can be expanded to a $\tau_{1}$-model $M^{+}$ from $K_{+}$, where: 

        \begin{itemize}
            \item $K_{+} = K_{\gk}^{+} = \{ N \colon N$ a $\tau_{1}$-structure omitting every  $p \in \cP$\}. 
        \end{itemize}

        \item[(B)] If $M^{+} \in K_{+}$ and $M^{+} \, {\rest} \, \tau \leq_{\gk} N$ then there is a model $N^{+} \in K_{+}$ expanding $N$ such that $M^{+} \subseteq N^{+}.$   Also, for $M, N \in K,$ we have $M \leq_{\gk} N$ \underline{iff} there are expansions $M^{+}, N^{+} \in K_{+}$ of $M, N$ respectively such that $M^{+} \subseteq N^{+}.$  

        \item[(C)] $(K_{+}, \subseteq)$ is an AEC with $\LST(K^{+}, \subseteq) = \LST(\gk).$ 

        \item[(D)] There is a set $\tau_{1}' \subseteq \tau_{1}$ of cardinality $\LST_{\gk}$ such that $A \subseteq M^{+} \in K_{+} \Rightarrow \rm{cl}_{\tau_{1}'}(A, M^{+}) \subseteq M^{+}.$

        \item[(E)] Some $\psi \in \bbL_{(2^{\lambda})^{+}, \aleph_{0}}$ defines $K_{\gk, 1}$ where $\lambda = \LST_{\gk} + \vert \tau_{\gk} \vert.$ 
    \end{enumerate}
\end{claim}

\begin{PROOF}{\ref{u11}}
    By \cite[1.7]{Sh:88r}.
\end{PROOF}

\begin{definition}\label{u14}\ 

    Assume $T$ is a theory in $\bbL_{\kappa, \aleph_{0}}(\tau_{T}),$ $\tau_{T}$ determined by $T$ (so $\vert T \vert \leq (\vert \tau_{T} \vert + \kappa)^{< \kappa}$) and recall $\mathcal{F}_{T}$ is the set of formulas $\varphi(\bar{x})$ such that $\varphi(\bar{x})$ is a  sub-formula of some sentence $\psi \in T.$ We define $\gk = \gk_{T}$ as follows:

   \begin{enumerate}
       \item[(A)] $K_{\gk}$ is the class of $\tau_{T}$-models of $T$ of cardinality $\geq \kappa + \vert T \vert$. 

       \item[(B)] $M \leq_{\gk} N$ \underline{iff}:

       \begin{enumerate}
           \item[(a)] $M, N \in K_{\gk},$

           \item[(b)] $M \subseteq N,$

           \item[(c)] $M \preceq_{\mathcal{F}} N$ i.e., if $\varphi(\bar{x}) \in \mathcal{F}_{T}$ (see below, so $\lg(\bar{x})$ is finite and $\bar{a} \in {}^{\lg(\bar{x})}M$) then $M \models \varphi[\bar{a}]$ \underline{iff} $N \models \varphi[\bar{a}].$ 
       \end{enumerate}
   \end{enumerate}
\end{definition}

\begin{claim}\label{u17}
    If $T$ is a theory in $\bbL_{\kappa, \aleph_{0}}(\tau_{T}),$ then: 

    \begin{enumerate}
        \item[(A)] $\gk_{T}$ is an AEC.

        \item[(B)] $\LST_{\gk_{T}} = \LST(\gk_{T}) \leq \vert T \vert + \kappa.$

        \item[(C)] If $T \subseteq \bbL_{\lambda^{+}, \aleph_{0}}(\tau_{T})$ then $\LST_{\gk} \leq \vert T \vert + \lambda.$ 

        \item[(D)] $\gk_{T}$ has no model of cardinality $< \vert \tau \vert + \kappa$ but for any $\tau(T)$-model $M$ of cardinality $\geq \vert T \vert + \kappa,$ $M \in K_{\gk_{T}} \Leftrightarrow M \models T.$
    \end{enumerate}
\end{claim}

\begin{PROOF}{\ref{u17}}
     Mainly, this holds by Fact \ref{X1.1} and Fact \ref{X1.2}, but see fully in the proof of Claim \ref{u26}.  
\end{PROOF}

\begin{conclusion}\label{u20}
    To prove our results for $T \subseteq \bbL_{\kappa, \aleph_{0}}$ it suffices to prove them for the AEC $\gk_{T}$ (see Definition \ref{u14}). 
\end{conclusion}

\begin{PROOF}{\ref{u20}}
    By Claim \ref{u17} just check the definitions and assumptions.  
\end{PROOF}

\begin{definition}\label{u23}
    We say the AEC  \emph{$\gk$ is $(\mu, \lambda, \kappa)$-representable} \underline{when} there are $(\tau_{1}, T_{1}, \Gamma)$ such that: 

    \begin{enumerate}
        \item[(a)] $\tau_{1} \supseteq \tau_{\gk}$ has cardinality $\leq \lambda,$

        \item[(b)] $T_{1} \subseteq \bbL(\tau_{1})$ is a first order logic universal theory, so $\vert T_{1} \vert \leq \lambda,$

        \item[(c)] $\Gamma$ is a set of $\leq \mu$ $\rm{qf}$-types in $\bbL(\tau_{1}),$ each of cardinality $< \kappa,$ 

        \item[(d)] $M \in K_{\gk}$ \underline{iff} $M$ is the $\tau_{\gk}$-reduct of some $M_{2} \in \EC(T_{1}, \Gamma),$ where $$ \EC(T_{1}, \Gamma) = \{ N  \colon N \text{ a } \tau_{1} \text{-model of } T_{1} \text{ omitting every } p(x) \in \Gamma \},$$

        \item[(e)] $M \leq_{\gk} N$ \underline{iff} for every $M_{1} \in \EC(T_{1}, \Gamma)$ expanding $M$, there is $N_{1} \in \EC(T_{1}, \Gamma)$ expanding $N$ and extending $M_{1}.$  
    \end{enumerate}
\end{definition}

\begin{claim}\label{u26}\ 

    (1) Let $\gk$ be an AEC. If $\lambda \geq \LST_{\gk} + \vert \tau_{\gk} \vert,$ then $\gk$ is $(2^{\lambda}, \lambda, \lambda^{+})$-representable.

    (2) If $T \subseteq \bbL_{\kappa, \aleph_{0}}$ is a theory then $\gk_{T}$ is $(\vert \tau \vert + \kappa, \vert \tau \vert + \kappa, \kappa)$-representable. If in addition $\vert T \vert < \kappa$ hence is $\subseteq \bbL_{\theta, \aleph_{0}}$ for some $\theta < \kappa$, then it is $(\vert \mathcal{F}_{T} \vert, \vert \mathcal{F}_{T} \vert, \theta)$-representable.
\end{claim}

\begin{PROOF}{\ref{u26}}\ 

    (1) By Claim \ref{u11}, that is, by \cite{Sh:88r} and classical theorem, see e.g. \cite[Ch.~VII]{Sh:f}. 

    (2) Just consider Definition \ref{u23} and the proof of Claim \ref{u17}. 
\end{PROOF}

\subsection{Indiscernibles and Ehrenfeucht-Mostowski structures}\label{1B}\ 

The basic results on generalized Ehrenfeucht-Mostowski models can be found in \cite{Sh:a} or \cite[VII]{Sh:c}. 

\begin{definition}\label{v0}\ 

    (1) We recall here some notation. Let {\bf I} be a class of models which we call the \emph{index models}. Denote the members of {\bf I} by $I, J\ldots,$ etc.
    
    (2) For $I\in$ {\bf I} we say that $\lng a_s \colon s\in I\rng$ is \emph{indiscernible in $M$} iff the $a_{s}$-s are pairwise distinct and for every $\bar s,\bar t \in {}^{ < \om} I$ realizing the same atomic type in $I$, $\bar a_{\bar s}$ and $\bar a_{\bar t}$ realize the same quatifier free type in $M$ (where $\bar {a}_{\lng s_0,\ldots,s_n\rng} = \langle a_{s_{0}}, \ldots,  a_{s_n} \rangle$). 
    
    (3) Assume $\tau \subseteq \tau'$ are vocabularies and $\Phi$ is a function with domain including $$\{ \tp_{\at}(\bar s,\emptyset,I) \colon \bar s\in {}^{<\om}I \text{ for some } I \in \bfI\}$$ 
    
    and if $\bar{s} \in {}^{n} I$ then $\Phi(\tp(\bar{s}, \emptyset, I))$ is a complete quantifier free $n$-type in $\bbL(\tau'),$ let $\tau_{\Phi} = \tau'$. Moreover, if $I\in$ {\bf I}, we let $\GEM'(I,\Phi)$ be an $\tau'$-model generated by
    $ \{ a_s : s \in I \} $ 
    such\footnote{Equivalently, we can use $\tp_{\mathrm{qf}},$ the quantifier free type.} that $\tp_{\at}(\bar a_{\bar s},\emptyset, M)=\Phi\big(\tp_{\at}(\bar s,\emptyset,I)\big);$ $\langle a_{s} \colon s \in I \rangle$ is called the skeleton. 
    
    (4) We say that $\Phi$ is \emph{proper for {\bf I}} if for every $I\in$ {\bf I}, $\GEM'(I,\Phi)$ is well-defined.

    (5) Let $\GEM(I,\Phi)$ be the $\tau$-reduct of $\GEM'(I, \Phi)$.

    Pedantically, we should write $\GEM_{\tau}(I, \Phi)$ but $\tau$ is constant. 
    
    (6) For the purposes of this paper we'll let {\bf I} be the class {\bf LO} of linear orders and $\Phi$ will be proper for {\bf LO} and then write $\EM$ (instead $\GEM$). For $I\in$ {\bf LO} we may abbreviate $\EM'(I,\Phi)$ by $\EM'(I)$ and $\EM(I,\Phi)$ by $\EM(I)$, when $\Phi$ is clear from the context. 
\end{definition}

We first deal with pairs $(T, \mathcal{F}).$ 

\begin{Claim}\label{v2} 
    If $T \subseteq \bbL_{\kappa, \aleph_{0}}(\tau)$ is a theory which has a model of cardinality $\geq \kappa$, then there are $\tau,$ $\Phi$ as in Definition \ref{v0} such that, for each linear order $I = (I, \le)$ there exists a Ehrenfeucht-Mostowski model $\EM(I, \Phi)$ is a model of $T.$
\end{Claim}

\begin{PROOF}{\ref{v2}}
    See Nadel \cite{Na85} and Dickmann \cite{Dic85} or \cite[VII, \S5]{Sh:c}.
\end{PROOF}

But now we use the AEC framework.

\begin{claim}\label{v8}\ 

    (0) If $\gk$ is an AEC then without loss of generality $\tau_{\gk}$ has cardinality $\leq 2^{\LS(\gk)}$. Fully we have $\tau_{\gk} / E_{\gk}$ has $\leq 2^{\LS(\gk)}$ equivalent classes when $E_{\gk} = \{ (R_{1}, R_{2}) \colon R_{1}, R_{2}$ are both predicates or both function symbols and are of the same arity and $M \in K_{\gk} \Rightarrow R_{1}^{M} = R_{2}^{M} \}$.  

    (1) Assume $\gk$ is an AEC, $\mu = 2^{\LST(\gk) + \vert \tau(\gk) \vert}$. If $\gk$ has a model of cardinality $\geq \beth_{\mu^{+}}$ (or just model of cardinality $\geq \beth_{\alpha}$ for every $\alpha < \mu^{+}$) \underline{then} there is $\Phi$ such that: 

    \begin{enumerate}
        \item[(a)] $\Phi$ is as in Definition \ref{v0}, 

        \item[(b)] $\tau_{\Phi} = \tau_{\gk, 1},$ where $\tau_{\gk, 1}$ is from Claim \ref{u11} or Definition \ref{u23}, 

        \item[(c)] $\EM(I) \in K_{\gk}$ has cardinality  
        $\LST_{\gk} + \vert I \vert,$

        \item[(d)] for $(\tau_{1}, T_{1}, \Gamma)$ as in Definition \ref{u23}, every model of the form $\EM'(I)$ is in $\EC(\Gamma, T_{1})$ and $\tau_{\Phi} = \tau_{1}.$
    \end{enumerate}

    (2) In particular, 

    \begin{enumerate}
        \item[(a)] $\EM'(I)$ is a $\tau_{1}$-model,

        \item[(b)] $\EM(I) = \EM'(I) \, {\rest} \, \tau$ belongs to $K,$

        \item[(c)] (follows) if $I \subseteq J$ then $\EM(I) \leq_{\gk} \EM(J),$ both models from $K$ of cardinality $\vert I \vert + \LST(\gk).$
    \end{enumerate}
\end{claim}

\begin{PROOF}{\ref{v8}}
    As in \cite[1.13]{Sh:88r}, \cite[Ch. VII]{Sh:c}. 
\end{PROOF}

\subsection{Limit ultra-powers, iterated ultra-powers and nice extensions}\label{1C}\ 

An important technique we shall use in studying the categoricity spectrum of a theory in $\bbL_{\k, \aleph_{0}}$ or  suitable AECs is the limit ultra-power.  It is convenient to record here the well-known definitions and  properties of limit and iterated ultra-powers (see  Chang and Keisler \cite{CK73}, Hodges-Shelah \cite{Sh:109}) and then to examine nice extensions of models.

\begin{Definition}\label{w2}
    Suppose that $M$ is an $\tau$-structure, $I$ is a non-empty set, $D$ is an ultra-filter on $I$ (but see Definition \ref{w5}(5)), and $G$ is a filter on $I \times I$.
    
    (1) For each $g \in {}^{I}|M|$, let
    
    \begin{enumerate}
        \item[(a)] $\eq(g) \coloneqq \{ \langle i, j \rangle \in I \times I \colon g (i) = g (j) \},$ and

        \item[(b)] $g / D \coloneqq \{ f \in {}^I |M| \colon g = f \,  \Mod \, D \}$ where, $$g = f \; \Mod \; D \text{ iff }\{ i \in I \colon g (i) = f (i) \} \in D.$$ 
    \end{enumerate}
    
    (2) Let $\underset{D / G} \Pi |M| \coloneqq \{ g / D \colon  g \in {}^I |M|$ and $\eq (g) \in G  \}$. Note that $\underset{D / G} \Pi |M|$ is a non-empty subset of $\Pi_D |M| = \{ g / D  \colon g \in {}^I|M| \}$ and is closed under the constants and functions of the ultra-power $\Pi_D M$ of $M$ modulo $D$. 
    
    (3) The limit ultra-power $\underset{D/G}  \Pi M$ of the $\tau$-structure $M$ (with respect to ($I, D, G$)) is the substructure of $\Pi_D M$ whose universe is the set $\underset{D / G} \Pi |M|$.  The canonical map $d$ from $M$ into $\underset{D / G} \Pi M$ is defined by $d (a) = \langle a_i : i \in I \rangle / D$, where $a_i = a$ for every $i \in I$.  
    
    (4) Note that the limit ultra-power $\underset {D / G}  \Pi M$ depends only on the equivalence relations which are
    in $G$, i.e. if {\bf E} is the set of all equivalence relations on $I$ and $G
    \cap$ {\bf E} $= G' \cap$ {\bf E}, where $G'$ is a filter on $I \times I$, then $\underset
    {D / G}  \Pi M = \underset{D / G'}  \Pi M$.
\end{Definition}

\begin{Definition}\label{w5}

    Assume, 
    
    \begin{enumerate}
        \item[(a)] $M$ be an $\tau$-structure, $\langle Y, < \rangle = \langle Y, <_{Y} \rangle$ a linear order,

        \item[(b)] for each $y \in Y$, let $D_y$ be an ultra-filter on a non-empty set $\langle Y, <_{Y} \rangle$, 

        \item[(c)] $\bar{I} = \langle I_{y} \colon y \in Y \rangle,$

        \item[(d)]  $ \bar{D} = \langle D_{y} \colon y \in Y \rangle$,

        \item[(e)] $I = \underset{y \in Y} \Pi I_y$.
    \end{enumerate}    

    Then, 

    (1) Let $E = \underset{y \in Y}{\Pi} D_y $ be the set of $ s \subseteq I$ such that there are $y_1 < \dots < y_n$ in $Y$ satisfying:
    
    \begin{enumerate}
        \item[$(\alpha)$] for all $i, j \in I$, if $i \uhr \{ y_1, \dots, y_n \} = j \uhr \{ y_1, \dots, y_n \}$ then $i \in s$ iff $j \in s,$
        
        \item[$(\beta)$] $\{ \langle i (y_1), \dots, i (y_n) \rangle: i \in s \} \in D_{y_1} \times \dots \times D_{y_n}.$
    \end{enumerate}
    
    (2) The iterated ultra-power $\prod_{\bar{D}} \vert M \vert$ or $\prod_E |M|$ of the set $|M|,$ noting $E$ is a filter on $I,$ is the set $\{ f /E \colon  f \in {}^{I} \vert M \vert$ and for some finite $Z_f \subseteq Y$ for all $i, j \in I$, if $i \uhr Z_f = j \uhr Z_f$, then $f (i) = f (j) \}$. 

    (2A) Note that $\langle Y, < \rangle, \bar{I}, \bar{D}, E, I$ and $E$ can be defined from $E$ and can be defined from $\bar{D},$ so we may write $\Pi_{E},$ $\Pi_{\bar{D}}$ above. 
    
    (3) The iterated ultra-power $\prod_E M$ of the $\tau$-structure $M$ with respect to $\langle D_y:  y \in Y \rangle$ is the $\tau$-structure whose universe is the set $\Pi_E |M|$; for each $n$-ary predicate symbol $R$ of $L, R^{\Pi_E M} (f_1 / E, \dots, f_n / E)$ iff $\{i \in I \colon R^M (f_1 (i), \dots, f_n (i )) \} \in E$; for each $n$-ary function symbol $F$ of $L ,\; F^{\Pi_E M} (f_1 / E, \dots, f_n / E) = \langle F^M (f_1 (i), \dots, f_n (i))  \colon i \in I \rangle / E$. 
    
    (4) The canonical map $d \colon  M \to \Pi_E M$ is defined as usual by: $$d (a) = \langle a  \colon  i \in H \rangle / E.$$

    (5) In  Definition \ref{w2}, we do not need ``$D$ is an ultra-filter on $I$'', just ``$D$ is a filter on $I$ such that, if $e \in G \cap \mathbf{E},$ then $D / e$ is an ultra-filter on $I / e$''.  

    (6) We say $\bfu$ is an iterated ultra-powers parameters when it consists of  $\langle Y, < \rangle$, $\bar{I} = \langle I_{y} \colon y \in Y \rangle$, $\bar{D} = \langle D_{y} \colon y \in Y \rangle$ and $I$ as in the beginning of Definition \ref{w5}, $E$ as in Definition \ref{w5}(1) and

    \begin{itemize}
        \item $G = \{ e \colon$ is an equivalence relation on $I$ such that, for some finite subset $Z$ of $Y$, we have $f, g \in I \wedge f \rest Z = g \rest Z \Rightarrow f \mathrel{e} g\}.$ 
    \end{itemize}

    (6A) So $\bfu = (Y_{\bfu}, \dots, \,)$ definable from $\bar{D}_{\bfu}$ and from $E_{\bfu}$ and we may write $\prod_{\bfu}.$
\end{Definition}

\begin{Remark}\label{w8}\ 

    (1) Every ultra-power is a limit ultra-power:  take $G =  \cP (I \times I)$ and note that $\Pi_D M = \underset{D/G} \Pi M$.
    
    (2) Every iterated ultra-power is a limit ultra-power, hence in Definition \ref{w5} we may write $\Op_{\bar{D}},$ $\Op_{E}$ or $\Op_{\bfu}.$ 
    
    [Why? let the iterated ultra-power be defined by $\langle Y, <\rangle$ and $\langle(I_y, D_y) \colon y\in Y\rangle$ (see Definition \ref{w5}). For $Z \in [Y]^{< \w}$, let  $A_Z = \{ (i,j) \in I \times I \colon  i \uhr Z = j \uhr Z \}$.  Note that
    $\{ A_Z  \colon Z \in [Y]^{< \w } \}$ has the finite intersection property and hence can be extended to a filter $G$ on $I \times I$.  Now for any model $M$ we have $\Pi_E M \cong \underset{D/G} \Pi M$  for every filter $D$  over $I$ extending $E$ under the map $f / E \to f / D$.]
\end{Remark}

\begin{Definition}\label{w11}\ 

    (1) We say that $(I, D, G)$ is \emph{suitable} when: 

    \begin{enumerate}
        \item[(a)] $D$ is an ultra-filter on a non-empty set $I$ (or just a filter, see Definition \ref{w5}(5)),

        \item[(b)] $G$ is a \underbar{suitable}, pedantically  a $D$-suitable filter on $I \times I$ or just a set of equivalence relations on $I$, which means:
    
        \begin{enumerate}
            \item[(i)]  if $e \in G$ and $e'$ is an equivalence relation on $I$ coarser than $e$, then $e' \in G,$
    
            \item[(ii)] $G$ is closed under finite intersections,
    
            \item[(iii)] $(I, D, G)$ is \emph{$\kappa$}-complete, which means that, if $e \in G$, then $D/e = \{ A \subseteq I/e: \underset{x \in A} \cup x \in D \}$ is a $\k$-complete ultra-filter on $I/e$ which, for simplicity, has cardinality $\kappa$.
        \end{enumerate}
    \end{enumerate}

    (2) A limit parameter $\bfu$ is \emph{suitable} when $(I_{\bfu}, E_{\bfu}, G_{\bfu})$ is. 
    
    (3) Suppose that $M$ is an $\tau$-structure and $(I, D, G)$ is suitable. \underline{Then} $\Op(M, I, D, G) = \Op_{I, D, G}(M)$ is the limit ultra-power $\underset{D/ \hat G} \Pi M$ where $\hat G$ is the filter on $I \times I$ generated by $G$. When clear from the context one abbreviates $\Op(M, I, D, G)$ by $\Op(M)$, pedantically $\Op$ stand for $\Op_{I, M, G}$ and one writes $f_{\Op} = f_{\Op, M}$ for the canonical map $d \colon M \to \Op(M);$ so we may write $f_{\Op}$ or $f_{\Op}^{M}$ instead $f_{\Op, M}$ when $M$ is clear from the context. 
\end{Definition}

Recall that,

\newtheorem{obscon}[theorem]{Observation / Convention}

\begin{obscon}\label{w14}\ 

    (1) For any $\tau$-structure $N, f_{\Op} = f_{\Op, N}$ is an elementary embedding of $N$ into $\Op(N)$ and in particular $f_{\Op}: N \underset{\gk} \rightarrow \Op(N)$.
    
    (2) Since $f_{\Op}$ is canonical, one very often identifies $N$ with the $\tau$-structure $\mathrm{rng}(f_{\Op})$ which is an $\gk$-elementary substructure of $\Op(N)$, and one writes $N \leq_{\gk} \Op(N)$.  In particular for any model $M \in K$  and $\Op, f_{\Op} \colon M  \to_{\gk} \Op(M)$ (briefly written $M \leq_{\gk} \Op(M)$) so that $\Op(M)$ is a model from $K$ too.
    
    (3) Remark that if $D$ is a $\k$-complete ultra-filter on $I$ and $G$ is a filter on $I \times I$, then $\Op(M, I, D, G)$ is well defined.
    
    (4) Suitable limit ultra-power means one using a suitable triple, for such $\Op$ in Observation/Convention \ref{w14}(2) we get a $\bbL_{\kappa, \aleph_{0}}$-elementary embedding.
\end{obscon} 

More information on  limit and iterated ultra-powers can be found in \cite{CK73} and \cite{Sh:109}. 

\begin{Observation}\label{w17}
    (1) Given $\kappa$-complete ultra-filters $D_1$ on $I_1, D_2$ on $I_2$  and suitable filters $G_1$ on $I_1 \times I_1, G_2$ on $I_2 \times I_2$ respectively, there exist a $\kappa$-complete ultra-filter $D$ on a set $I$ and a filter $G$ on $I \times I$ such that: $$ \Op (M, I, D, G) = \Op(\Op(M, I_1, D_1, G_1), I_2, D_2, G_2) $$ and $(D, G, I)$ is $\kappa$-complete. 
    
    (2) Also iterated ultra-power (along any linear order) with each iterand being ultra-power by $\kappa$-complete ultra-filter, gives a suitable triple (in fact, even iteration of suitable limit ultra-powers is a suitable ultra-power).
\end{Observation}

\begin{Definition}\label{w20}
    Suppose that $K$ is a class of $\tau$-structures and $< \, = \, <_{\mathcal{K}}$ is a binary relation on $K$ (usually $(K, <) = (K_{\gk}, <_{\gk})$).  For $M, N \in K$, write $f \colon  M \leq_{\mathcal{K}}^{\mathrm{nice}} N$ to mean (if $<$ is clear from the context we may write $f \colon M \underset{\rm{nice}}{\to} N$ and, if $f = \id_{M}$ we may write $M \underset{\rm{nice}}{\leq} N$):

    \begin{enumerate}
        \item[(a)] $f$ is an isomorphism from $M$ into $N$ and $\mathrm{rng}(f) < N.$ Which means $f(M) < N,$ where $f(M)$ is the model $M'$ with universe $\mathrm{rng}(f)$ such that $f$ is an isomorphism from $M$ into $M',$ 

        \item[(b)] for some\footnote{We could use here and Theorem \ref{X2.5} suitable tuples $(I, D, G)$. However, then we have to add to the definition of ``$\gk$ is $(I, D, G)$-compact'' a clause saying: 
        
        \begin{enumerate}
            \item[$(\ast)$] if $M \in K_{\gk},$ $e_{1} \supseteq e_{2}$ are from $G$ and $M_{\ell} = \prod_{D / G} M \rest \{ f \in {}^{I}M \colon \eq(f) \supseteq e_{\ell} \}$ for $\ell = 1, 2$, \underline{then} $M_{1} \leq_{\gk} M_{2}$. 
        \end{enumerate}
        
        In \cite{Sh:362} this issue does not arise.} ultra-limit parameter $\bfu = (Y, <_{Y}, \bar{I}, \bar{D}, I, E, G)$, so $G$ is a suitable set of equivalence relations on $I$ (so Definition \ref{w11} clause (i), (ii), (iii) holds), and an isomorphism $g$ from $N$ into $\Op(M, I, E, G)$ such that $\mathrm{rng}(g) < \Op(M, I, E, G)$ and $g f = f_{\Op}$, where $f_{\Op}$ is the canonical embedding of $M$ into $\Op(M, I, E, G)$. Then $f$ is called a $<$-nice embedding of $M$ into $N$. Of course, one writes $f \colon M \underset{\text{nice}}  \rightarrow N$ and says that $f$ is a nice embedding of $M$ into $N$ when $<$ is clear from the context.
    \end{enumerate}
\end{Definition}

\begin{Example}\label{X1.9.1}\
    Consider $T, \Cal F$ and $\Cal K = \langle K, < \rangle = (K, <_{\mathcal{K}})$ as set up in Example \ref{X1.3B}.  In this case $f \colon M \underset{\rm{nice}} {\rightarrow} N$ holds iff $f \colon M \underset {\Cal F}  \rightarrow N$ and for some suitable ultra-limit parameter $\bfu$ and some $g \colon N \underset {\Cal F}  \rightarrow \Op_{\bfu}(M) $ we have $ g f = f_{\Op}$.
\end{Example}

Abusing notation one may writes $M \underset{\rm{nice}} \rightarrow N$ to mean that there are $f, g$ and $\Op$ such that $f \colon M  \underset{\rm{nice}}{\rightarrow} N$ using $g$ and $\Op.$  IF NOT SAID OTHERWISE, $<$ is $<_{\gk}$. We may also write $M{\underset{\rm{nice}}  \le } N$, and for linear orders we use $I \underset{\text{nice}} \subseteq J$.

\begin{Example}\label{X1.9.2}
    Let {\bf LO} be the class of linear orders and let $(I, \le_I) < (J, \le_J)$ mean that $(I, \le_I) \subseteq (J, \le_J)$, i.e. $(I, \le_I)$ is a suborder of $(J, \le_J)$. If $f: (I, \le_I)  \underset{\rm{nice}}{\rightarrow} (J, \le_J)$, then identifying isomorphic orders, one has $(I, \le_I) \subseteq (J, \le_J) \subseteq \Op(I, \le_I)$ and we 
    may write $(I, \leq_{I}) \underset{\rm{nice}}{\subseteq} (I, \leq_{J}).$
\end{Example}

\begin{Observation}\label{X1.10}
    Assume that $\gk$ is  as in \ref{u14}.  Suppose further  $M  \underset{\mathrm{nice}}{\leq} N$ and $M \subseteq M' \leq_{\gk} N$ where $M, M', N \in K$.  Then $M \underset{\mathrm{nice}}{\leq} M'$. 
\end{Observation}

\begin{PROOF}{\ref{X1.10}}
    For some $f, g$ and $\Op$, $f \colon  M \underset{\gk} \rightarrow N$, $g \colon N \underset{\gk} \rightarrow \Op(M)$ and $g f = f_{\Op}$.  Now $g \colon M' \underset{\gk} \rightarrow \Op(M)$ (since $M' \leq_{\gk} N)$ and $gf = f_{\Op}$ so that $M \underset{\rm{nice}}  \leq 
    M'$.
\end{PROOF}

\begin{observation}\label{X1.11}
    Suppose that $\delta$ is any ordinal, $\lng M_i \colon i\le\del\rng$ is a continuous increasing chain and for each $i<\del$, $M_i\underset{\text{nice}}  \leq M_{i+1}$. Then for every $i<\del$, $M_i\underset{\rm{nice}}  \leq M_\del$.
\end{observation}

\begin{PROOF}{\ref{X1.11}}
    Like the proof of Remark \ref{w8}(2). For each $i<\del$,  there is a $\bfu_{i}$ as in Definition \ref{w11} which witnesses $M_i\underset{\text{nice}} \le M_{i+1}$ and
    and let 
    $Y_{i} = Y_{\bfu_{i}}$ for $i < \delta.$ Without loss of generality, $\langle Y_{i} \colon i < \delta \rangle$ are pairwise disjoint. We define $\bfu$ by: 

    \begin{enumerate}
        \item[(a)] $Y = \bigcup \{ Y_{i} \colon i < \delta \}$,

        \item[(b)] $s <_{y} t$ iff $\bigvee_{i < \delta} s <_{i} t$ or $s \in Y_{\bfu_{i}} \wedge t \in Y_{\bfu_{j}} \wedge i<j,$

        \item[(c)] $D_{j} = D_{\bfu_{i}, s}$ when $s \in Y_{\bfu_{i}}$ for $i < \delta$.
    \end{enumerate}

    This is enough and the rest should be clear. 
\end{PROOF}

\begin{Claim}\label{X1.12}
    For every model $M$ of cardinality $\geq \kappa$ and $\lambda \geq \kappa + \LST_{\gk} + \| M\|$ there is $N$ such that
    $M \underset{\text{nice}}  {\leq} N$, $M\neq N$ and $\|N\|=\lambda$.
\end{Claim}

\begin{PROOF}{\ref{X1.12}} 
    As $\gk$ is $\bfD$-compact, by Assumption \ref{X1.-2}(2) no $M \in K_{\geq \kappa}$ is $\leq_{\gk}$-maximal, so by Definition \ref{u2} we are done. 
\end{PROOF}

\begin{definition}\label{y1}\ 

    (1) Assume $D$ is an ultra-filter on $\kappa.$ For an AEC $\gk = (K_{\gk}, \leq_{\gk})$ we say $\gk$ is \emph{$D$-compact} \ 
    \underline{when}: 

    \begin{enumerate}
        \item[(a)] if $M \in K_{\gk}$ then the ultra-power $M^{\kappa} / D$ belongs to $K_{\gk},$ 

        \item[(b)] moreover, the canonical embedding of $M$ into $M^{\kappa} / D$ is a $\leq_{\gk}$-embedding, 

        \item[(c)] if $M \leq_{\gk} N$ then the canonical embedding of $M^{\kappa} / D$ into $N^{\kappa} / D$ is a $\leq_{\gk}$-embedding,

        \item[(d)] $\gk$ has a model of cardinality $\geq \kappa$ (or at least of cardinality $\geq \theta$ where $D$ is not $\theta$-complete). 
    \end{enumerate}        
        
    (2) If $\bfu = (Y, \bar{I}, \bar{D}, I, E, G)$ is as in Definition \ref{w5}, then for an AEC $\gk$ we say $\gk$ is $\bfu$-compact and $E$-compact  \underline{when}: 

    \begin{enumerate}
        \item[(a)] if $M \in K_{\gk}$ and $\prod_{E}M \in K_{\gk},$

        \item[(b)] moreover, the canonical embedding of $M$ into $\prod_{E}M$ is a $\leq_{\gk}$-embedding,

        \item[(c)] if $M \leq_{\gk} N$ then the canonical embedding of $\prod_{E} M$ into $\prod_{E} N$ is a $\leq_{\gk}$-embedding. 
    \end{enumerate}
\end{definition}

\begin{claim}\label{y3}

    Assume $D$ is a non-principal $\kappa$-complete ultra-filter (usually on $\kappa$). 

    (1) If $\gk$ is an AEC and $\vert \tau_{\gk} \vert + \LST(\gk) < \kappa$ \underline{then} $\gk$ is $\bfD$-compact.

    (2) If $\gk$ is $(\mu, \lambda, \kappa)$-representable, \underline{then} $\gk$ is $D$-compact.

    (3) Also the claim on $\Op$ generalizes, that is, if $ \langle Y, < \rangle, \bar{I}, \bar{D}, E, I$ is as in Definition \ref{w5} and $\gk_{s}$ is $D_{s}$-compact for every $s \in Y$ then in (1) and (2), $\gk$ is $E$-compact. 

    (4) So if there is one ultra-filter $D$ on $\kappa$ which is normal or just non-principal $\kappa$-complete ultra-filter on $\kappa$, \underline{then} for every linear order $\langle Y, < \rangle$ then we can find $\bar{I}, \bar{D}, E, I$ such that they together are as in \ref{y3}(3)  
\end{claim}

\begin{PROOF}{\ref{y3}}\ 

    (1) By \ref{u26} and part (c).

    (2), (3), (4) Easy. 
\end{PROOF}

\begin{claim}\label{b0}
    Assume $D$ is a non-principal $\kappa$-complete ultra-filter on $\kappa$ and $\gk_{1}$ is a $D$-compact AEC, $\chi \geq \LS_{\gk_{1}}$ and let $\gk_{2} = (\gk_{1})_{[ \chi, \infty )},$ see Definition \ref{u5}(3). 

    (1) If $\chi \geq \kappa$ and $\gk_{1}$ is $D$-compact then $\gk_{2}$ is $D$-compact. 

    (2)  If $\lambda \geq \chi,$ then $\gk_{1}$ is categorical in $\lambda$ iff $\gk_{1}$ is categorical in $\lambda.$

    (3) If $\lambda \geq \chi,$ then $(\gk_{2})_{[ \chi, < \lambda)}$ has amalgamation iff $(\gk_{1})_{[ \chi, < \lambda )}$ has amalgamation. 
\end{claim}

\begin{PROOF}{\ref{b0}}
    Straightforward. 
\end{PROOF}

\begin{remark}\label{b1}\ 

    (1) Claim \ref{y3} justifies the assumption $\LS_{\gk} \geq \chi$ in \ref{b2} below (e.g. to proof \ref{u8}). 

    (2) Usually $\lambda$ denotes a power in which $\gk$ is categorical. 
\end{remark}

For the rest of this work,

\begin{hypothesis}\label{b2}
    Assume $\chi \geq \kappa.$ 

    (1) $\gk$ is a $\bfD$-compact  AEC with $\LS_{\gk} = \chi$, no $M \in K_{\gk}$ has cardinality $< \chi$, $\bfD$ a $\kappa$-complete non-principal ultra-filter on $\kappa$, $K = K_{\gk}$ and similarly for any $\langle Y, < \rangle, \bar{I}, \bar{D}, I$ or $E$ derived from $\bfD$ as in Definition \ref{w5}. 
    
    (2) $\Phi,$ $\bfa = \langle a_{s} \colon s \in I \rangle$ \ are as in Definition \ref{v0} for $\gk$ with $\tau_{\Phi}$ of cardinality $\leq \chi$, hence $\lambda \geq \chi  \Rightarrow (K_{\gk})_{\lambda} \neq \emptyset$. 

\end{hypothesis}

\newpage

\section{The amalgamation property for regular categoricity}\label{2}

The main aim of this section is to show that if $K$ is categorical in the regular cardinal $\lambda > \LST_{\gk}$, then $\gk_{< \lambda} = \langle K_{< \lambda}, \leq_{\gk} \rangle$ has the amalgamation property (AP) (Definition \ref{X1.3} (1)). Categoricity is not presumed if not required.

Recall Hypothesis \ref{b2} is assumed. 

\begin{Lemma}\label{X2.1}
    Suppose that $\chi \le \mu \le \lambda,  M, M_1, M_2 \in K_{\le \mu},  f_1 \colon M \underset{\text{nice}}  \rightarrow M_1, f_2: M \underset{\Cal \gk} \rightarrow M_2$.  \underline{Then} there is an amalgam $N \in K_{\le \mu}$ of $M_1, M_2$ over $M$ with respect to $f_1, f_2$.
    
    Moreover, there are $N$ and $g_\ell \colon M_\ell \underset{\gk} \rightarrow N$ for $\ell = 1, 2$ such that $g_1f_1 = g_2 f_2$ hence $\mathrm{rng}(g_2 f_{2}) = \mathrm{rng}(g_1f_1)$ and $g_{1} \colon M_{1} \underset{\rm{nice}}{\to} N.$
\end{Lemma}

\begin{PROOF}{\ref{X2.1}}
    There are $g$ and $\Op$ such that $g \colon M_1 \underset{\gk} \rightarrow \Op(M),  g f_1 = f_{\Op, M}$. Now, $f_{2}$ induces an $\leq_{\gk}$-elementary embedding $f^\ast_2$ of $\Op(M)$ into $\Op(M_2)$ such that $f^\ast_2 f_{\Op}^{M} = f_{\Op}^{M_{2}} f_2$.  Let $g_1 = f^\ast_2 g$ and $g_2 = f_{\Op, M_{2}}$. By Fact \ref{X1.2} one finds $N \in K_{\le_\mu}$ such that $\mathrm{rng}(g_1) \cup \mathrm{rng}(g_2) \subseteq N \leq_{\gk} \Op(M_2)$.  Now $N$ is an amalgam of $M_1, M_2$ over $M$ with respect to $f_1, f_2$ since $g_1 f_1 = f^\ast_2 g f_1 = f^\ast_2 f_{0_p, M} =
    f_{\Op, M_{2}} f_2 = g_2 f_2$. The last phrase in the lemma is easy by properties of $\Op$.
\end{PROOF}

\begin{Lemma}\label{X2.2}
    Suppose that $M \in K_{\le \mu}$ is a $\mu$-c.a.b., $\chi \le \mu < \lambda$.  Then $N \in K_{< \lambda}$ is a $\| N \|$-c.a.b. whenever $f \colon M \underset{\text{nice}}  \rightarrow N$.
\end{Lemma}

\begin{PROOF}{\ref{X2.2}}
    By the assumption, there is $g \colon N \underset{\gk} \rightarrow \Op(M) $ such that $ g f = f_{\Op, M}.$ Recall  $M$ is a $\mu$- c.a.b., so for some $M_i \in K_{\le \mu}$ and $f_i \colon M \underset{\gk} \rightarrow M_i $ (for $i = 1,2)$ there is no amalgam of $M_1, M_2$ over $M$ w.r.t. $f_1, f_2$. Let $f^\ast_i$ be the $\leq_{\gk}$-elementary embedding from $\Op(M)$ into $\Op(M_i)$ induced by $f_i$ (note that $f^\ast_i f_{\Op, M} = f_{\Op, M_{i}} f_i,  i = 1,2)$. Choose $N_i$ of power $|| N ||$ such that $M_i \cup \mathrm{rng}(f^\ast_i g) \subseteq N_i \leq_{\gk} \Op(M_i)$.  Note that $f^{\ast}_{i} g \colon N \underset{\gk} \rightarrow N_i$. It suffices to show that there  is no amalgam of $N_1, N_2$ over $N$ w.r.t. $f^\ast_1 g, f^\ast_2 g$.
    
    Well, suppose that one could find an amalgam $N^\ast$ and $h_i \colon N_i \underset{\gk} \rightarrow N^\ast, i = 1,2$, with $h_1 (f^\ast_1 g) = h_2 (f^\ast_2 g)$.  Using Fact \ref{X1.2} choose $M^\ast, || M^\ast|| \le \mu,  M^\ast \leq_{\gk} N^\ast, \mathrm{rng}(h_1 f_{\Op} \uhr M_1) \cup \mathrm{rng}(h_2 f_{\Op} \uhr M_2) \subseteq | M^\ast|$.
    Set $g_i = h_i f_{\Op} \uhr M_i$, for $i = 1,2$, and note that: 
    \begin{equation*}
        \begin{split}
            g_1 f_1 &= h_1 f_{\Op} f_1 = h_1 f^\ast_1 f_{\Op} = h_1 f^\ast_1 g f = h_2 f^\ast_2 g f = h_2 f^\ast_2 f_{\Op}\\
            & = h_2 f_{\Op} f_2 = g_2 f_2.
        \end{split}
    \end{equation*}
    
    In other words, $M^\ast$ is an amalgam of $M_1, M_2$ over $M$ w.r.t. $f_1, f_2$-contradiction.  It follows that $N$ is a $\| N \|$-c.a.b.
\end{PROOF}

\begin{Corollary}\label{X2.3}
    Suppose that $\mu, \lambda$ satisfy $\chi \le \mu < \lambda$.  If $M \in K_\mu$ is a $\mu$-c.a.b., then there exists $M^\ast \in K_\lambda$ such that:
    
    \begin{itemize}
        \item[$(\ast)$]
             $M \leq_{\gk} M^*$ and for every $M' \in K_{< \lambda}$, if $M \leq_{\gk} M' \leq_{\gk} M^\ast$, then $M'$ is a $|| M'||$-c.a.b.
    \end{itemize}
\end{Corollary}

\begin{PROOF}{\ref{X2.3}}
    As $\| M \| \ge \k$, for some appropriate $\Op$ one has $|| \Op(M) || \ge \lambda$, and by Fact \ref{X1.2} one finds $M^\ast \in K_\lambda$ such that $M \subseteq M^\ast \leq_{\gk} \Op(M),$ hence $M \leq_{\gk} M^{\ast}$. Let us check that $M^\ast$ works in $(\ast)$.  Take $M' \in K_{< \lambda},  M \leq_{\gk} M' \leq_{\gk} M^\ast$; so $M \underset{\text{nice}}  \leq M'$ since $M^\ast \leq_{\gk} \Op(M),$ see Observation \ref{X1.10};  hence by Lemma \ref{X2.2}, $M'$ is a $|| M'||$-c.a.b.
\end{PROOF}

\begin{Theorem}\label{X2.4}
    Suppose that $\gk$ is $\lambda$-categorical, $\lambda = \cf(\lambda) > \chi$.  If $K_{< \lambda}$ fails AP, \underline{then} there is $N^\ast \in K_\lambda$ such that for some continuous increasing $\leq_{\gk}$-chain $\langle N_i \in K_{< \lambda}: i < \lambda \rangle$ of models,
     
    (1) $N^\ast = \underset{i < \lambda} \cup N_i,$
    
    (2) for every $i < \lambda, N_i \underset{\text{nice}}  \nleq N_{i+1}$
    (and so $N_i \underset{\text{nice}}  \nleq N^\ast$).
\end{Theorem}

\begin{PROOF}{\ref{X2.4}} 
    By an assumption $\gk_{< \lambda}$ fails AP, so for some $\mu \in [\chi, \lambda)$ and $M \in K_{\leq \mu},  M$ is a $\mu$-c.a.b. recalling Definition \ref{X1.4}. By Lemma \ref{X2.2} and Claim \ref{X1.12} without loss of generality $M\in K_{\mu}$. Choose by induction a continuous strictly increasing
    $\leq_{\gk}$-chain $\langle N_i \in {\text{K}}_{< \lambda} \colon i < \lambda \rangle$ as follows:
    
    $N_0 = M$; at a limit ordinal $i$, take the union; at a successor ordinal $i = j +1$, if there is $N \in K_{< \lambda}$ such that $N_j \leq_{\gk} N$ and $N_j \underset{\text{nice}}  \nleq N$ (so necessarily $N_{j} <_{\gk} N$), choose $N_i = N$, otherwise choose for $N_i$ any non-trivial $\leq_{\gk}$-elementary extension of $N_j$ of power
    less than $\lambda$. Next, we prove:
    

    \begin{enumerate}
        \item[$\boxplus$]  $(\exists j_0 < \lambda)(\forall j \in (j_0, \lambda)) (N_j$ is a $|| N_j ||$-c.a.b.).
    \end{enumerate}

    Why $\boxplus$ holds? Suppose not.  So one has a strictly increasing sequence $\langle j_i \colon  i < \lambda \rangle$ such that for each $i < \lambda, N_{j_i}$ is not a $\| N_{j_i} \|$-c.a.b. Let $N_\ast = \underset{i < \lambda} \cup N_{j_i}$. So $\| N_\ast \| = \lambda$. Applying \ref{X2.3} one can find $M^\ast \in K_\lambda$ and $M \in K_{< \lambda}$ such that $M \leq_{\gk} M^{\ast}$ and whenever $M' \in K_{< \lambda}$ and $M \leq_{\gk} M' \leq_{\gk} M^\ast$, then $M'$ is a $\| M'\|$-c.a.b.
    
    Since $\gk$ is $\lambda$-categorical, there is an isomorphism $g$ of $N_\ast$ onto $M^\ast$.
    Let $N = g^{-1} (M)$ and $M_i = g^{-1} (N_i)$ for $i < \lambda$. Now, $|| N || = \mu < \cf(\lambda) = \lambda$, so there is $i_0 < \lambda$ such that $N \subseteq N_{j_{i_0}}$.
    
    In fact $N_{j_{i_0}}$ is a $|| N_{j_{i_0}} ||$-c.a.b.
    [Otherwise, consider $M_{j_{i_0}}$. Since $M \leq_{\gk} M_{j_{i_0}} \leq_{\gk} M^\ast$ and $\| M_{j_{i_0}} \| < \lambda, M_{j_{i_0}}$ is a $\| M_{j_{i_0}} \|$-c.a.b., so there are $f_\ell\colon M_{j_{i_0}} \underset{\Cal F}  \rightarrow M'_\ell, (\ell = 1,2)$, with no amalgam of $M'_1, M'_2$ over $M_{j_{i_0}}$ w.r.t. $f_1, f_2$. If $N_{j_{i_0}}$ is {\it not} a $|| N_{j_{i_0}} ||$-c.a.b., then one can find an amalgam $N^+ \in K_{\le || N_{j_{i_0}} ||}$ of $M'_1, M'_2$ over $N_{j_{i_0}}$ w.r.t. $f_1 g, f_2 g$ such that $h_\ell \colon M'_\ell \underset{\gk} \rightarrow N^+$ and $h_1 (f_1 g) = h_2 (f_2 g)$;
    so $h_1 f_1 = h_2 f_2$ and $N^+$ is thus an amalgam of $M'_1, M'_2$ over $M_{j_{i_0}}$ w.r.t. $f_1, f_2,  \| N^+ \| \le \| N_{j_{i_0}} \| = || M_{j_{i_0}} ||$-contradiction.]  This contradicts the choice of $N_{j_{i_0}}$.  So the statement $\boxplus$ is correct.
    
    It follows that for each $j \in (j_0, \lambda)$ there are $N^1_j, N^2_j$ in ${\text{K}}_{< \lambda}$ and $f_\ell \colon N_j \to_{\gk} N^\ell_j$ such that no amalgam of $N^1_j, N^2_j$ over $N_j$ w.r.t. $f_1, f_2$ exists.  By Lemma \ref{X2.1} for both $\ell \in \{ 1,2 \}, N_j \underset{\rm{nice}}{\nleq} N^\ell_{j+1}$.  So by the inductive choice of $\langle N_{j +1} \colon j < \lambda \rangle, \forall j \in (j_0, \lambda) (N_j \underset{\text{nice}}  \nleq N_{j+1})$.  Taking $N^\ast = \underset {j_0 < j < \lambda}  \cup N_j$, one completes the proof  (of course for $j_0 < j < \lambda, N_j \underset{\text{nice}} \nleq N^\ast$: if
    $N_j \underset{\mathrm{nice}}{\leq} N^\ast \leq_{\gk} \Op (N_j)$, then by Observation \ref{X1.10} $N_j \underset{\text{nice}}  \leq N_{j+1}$-contradiction).
\end{PROOF}

\begin{Theorem}\label{X2.5}
    Suppose that $I = (I, <_I), $  $J = (J, <_J)$ are linear orders and $I$ is a suborder of $J$. Let $\EM'(I, \Phi)$ be as in Definition \ref{v0}, so let $\langle a_{s}^{1} \colon s \in I \rangle$ be a skeleton of $M_{1}' = \EM'(I) = \EM'(I, \Phi)$, a $\tau_{\Phi}$-model, $\langle a_{s}^{1} \colon s \in I \rangle$ is an indiscernible sequence in $\EM(I)$ which generates it. Similarly, $M_{2}' = \EM'(J, \Phi)$, $\langle a_{s} \in s \in J \rangle$ and as standard, we assume $M_{1}' \subseteq M_{2}'$, $ s \in I \Rightarrow a_{s}^{1} = a_{s}$, let $M_{\ell} = \EM(I) = M_{\ell}' \rest \tau_{\gk}$. If $(I, <_I) \underset{\text{nice}}  \subseteq (J, <_J)$. Then $\EM(I)  \underset{\text{nice}}  \leq \EM(J)$.
\end{Theorem}

\begin{PROOF}{\ref{X2.5}}
    So there is a suitable ultra-limit\footnote{We write $I_{\bfu} = I(\bfu)$ to distinguish it from $(I, <_{I})$.} parameter $\bfu = (Y, <_{Y}, \bar{I}, \bar{D}, I_{\bfu}, E, G)$ witnessing $(I, <_{I}) \underset{\mathrm{nice}}{\subseteq} (J, <_{J})$, that is, we have $(I, \leq_{I}) \subseteq \Op_{I, D, G}((I, <_{I}))$ and $(J, <)$ is isomorphic over $(I, \leq)$ to some $(J', <)$ such that $(I, <_{I}) \subseteq (J', <) \subseteq \Op_{I, D, G}((I, <))$ and let $\pi$ be such isomorphism. 
    So for each $t \in J$, there exists $f_t \in {}^{I(\bfu)} I$ such that $\pi(t) = f_t / D$.  Note that if $t \in I$, then $f_t / D = f_{\Op} (t)$ so that without loss of generality for all $i \in I_{\bfu}, f_t (i) = t$.  Define a map $h$ from $\EM(J)$ into $\Op(\EM(I))$ as follows.  An element of $\EM(J)$ has the form $$ \sigma^{\EM' (J) } (a_{t_1}, \dots, a_{t_n}),$$ 
    
    where $t_1, \dots, t_n \in J$, $\sigma$ an $\tau_{\Phi}$-term. Define, for $t \in J, g_t \in {}^I \EM(I)$
    by $g_t (i) = a_{f_t (i)}$. 
    
    Note that $f_t (i) \in I$, so that $a_{f_t (i) } \in \EM(I)$ and so $g_t / D \in \Op(\EM(I) )$.
    Let $h ( \sigma^{\EM' (J) } (a_{t_1}, \dots, a_{t_n} ) )
    = \sigma^{\Op(\EM' (I) ) } (g_{ t_1} / D, \dots, g_{t_n} / D)$ which
    is an element in $\Op(\EM(I))$.  The reader is invited to check that
    $h$ is an $\leq_{\gk}$-elementary embedding of $\EM(J)$ into $\Op(\EM(I))$, and consequently $\EM(I)  \leq_{\gk} \EM(J)$, but we elaborate. Prove by induction on $n < \omega$ that:  

    \begin{enumerate}
        \item[$\oplus$] if $\bar{s} = \langle s_{i} \colon i < n \rangle$ is $<_{Y}$-increasing then let $m \leq n$ and $N_{\bar{s}} \coloneqq M_{2} \rest \{ f \in {}^{H}M \colon$ $\eq(f)$ is refined by $\eq_{\bar{s} \rest m} = \{ (h_{1}, h_{2}) \colon h_{1}, h_{2} \in \prod_{s \in I} I_{s}$ and $\ell < m \Rightarrow h_{1}(s_{\ell}) \geq h_{2}(s_{\ell}) \} \},$ 
    \end{enumerate}

    \begin{enumerate}
        \item[$\boxplus$] for $\bar{s} = \langle s_{\ell} \colon \ell < n \rangle$ as above, $N_{\bar{s} \rest m} \leq_{\gk} N_{\bar{s}}.$ 
    \end{enumerate}

    [Why? Prove by induction on $n$ that it suffices to conclude that $m = n -1$ and now read the Definition.]

    \begin{enumerate}
        \item[$\boxplus$] if $\bar{s}$ is as above and $\bar{t}$ is a sub-sequence of $\bar{s}$ then $N_{\bar{t}} \leq_{\gk} N_{\bar{s}}$.
    \end{enumerate}
    
    Why? By \underline{Ax. V} of AEC (see Definition \ref{u2}): The rest should be clear. 
    
    Finally note that if $b = \sigma^{\EM' (I)} (a_{t_1}, \dots, a_{t_n}) \in \EM(I), t_1, \dots, t_n \in I$, then $h(a) = \sigma^{\Op(\EM'(I))} (g_{ t_1} / D, \dots, g_{t_n} / D) = \sigma^{\Op(\EM' (I))} ( \langle a_{f_{t_{i}(i)}}\colon i < \mu \rangle / D, \dots, \langle a_{f_{t_{n}(i)}} \colon i < \mu \rangle / D)  = $  $f_{\Op} (\sigma^{\EM' (I)} (a_{t_1}, \dots, a_{t_n})) = f_{\Op} (b)$.
    Thus $\EM(I) \underset{\text{nice}} \leq \EM(J)$.
\end{PROOF}

\begin{Criterion}\label{X2.6}
    Suppose that $(I, <)$ is a suborder of the linear order $(J, <)$. We have $(I, \leq) \underset{\rm{nice}}{\subseteq} (J, <)$ when:

    \begin{enumerate}
        \item[$(\ast)$] for every $t \in J \setminus I,$
    \end{enumerate}
    
    \begin{enumerate}
        \item[$(\aleph)$]  $\cf((I, <) \rest \{ s \in I\colon (J, <) \models s < t \}) = \kappa,$
    \end{enumerate}
    
    or
    
    \begin{enumerate}    
        
        \item[($\beth$)] $\cf((I, <)^{*} \rest \{  s \in I\colon (J, <)^{*} \models s <^{*} t\}) = \kappa.$
    \end{enumerate}
\end{Criterion}

\begin{notation}
      $(I, <)^\ast$ is the (reverse) linear order $(I^\ast, <^\ast)$ where $I^\ast = I$ and $(I^\ast, <^\ast) \vDash s <^\ast t$ iff $(I, <) \vDash t < s$.
\end{notation}

\begin{PROOF}{\ref{X2.6}}
    Let us list some general facts which facilitate the proof.
    
    Fact (A): Let $\underline{\k}$ denote the linear order $(\k,<)$ where $<$ is the usual order $\in \uhr \k \times\k$. If $J_1=\unl\k+J_0$, then $\unl\k {\underset{\text{nice}}  \subseteq} J_1$
    ($+$ is the addition of linear orders in which all elements in the first order
    precede those in the second).
    
    Fact (B): If $\unl\k \subseteq (I,<)$, $\unl\k$ is unbounded in $(I,<)$ and
    $J_1=I+J_0$, then $I\subnice J_1$.
    
    Fact (C): If $I\subnice J$, then $I+J_1\subnice J+J_1$.
    
    Fact (D): $I\subnice J$ iff $(J<)^*\subnice (I,<)^*$.
    
    Fact (E): If $\lng I_{\al} \colon\al \le \d\rng$ is a continuous increasing
    sequence of linear
    orders and for $\al<\d$, $I_{\al}\subnice I_{\al +1}$, then
    $I_\al\subnice I_{\d}$.
    
    Now using these facts, let us prove the criterion. Define an equivalence
    relation $E$
    on $J\setminus I$ as follows: $tEs$ iff $t$ and $s$ define the
    same Dedekind cut
    in $(I,<)$. Let $\{ t_{\al} \colon  \al<\d \} $ be a set of representatives of
    the $E$-equivalence
    classes. For each $\b \le \d$, define
    $$I_{\b}=J\uhr \left \{ t \colon t\in I \vee {\underset{{\al<\b}} \bigvee} tEt_{\al} \right \}$$
    so $I_0=I$, $I_{\d}=J$ and $\lng I_{\al} \colon\al \le \d\rng$ is a
    continuous increasing sequence of linear orders.
    By Fact (E), to show that $I \subnice J$,
    it suffices to show that $I_{\al}\subnice I_{\al+1}$ for each $\al <\d$.
    
    Fix $\al <\d$. Now $t_{\al}$ belongs to $J\setminus I$, so by
    $(*)$, $(\a)$ or $(\beth)$
    holds. By Fact (D), it is enough to treat the case $(\a)$. So without loss of generality
    $\Op \big( (I,<)\uhr \{ s\in I \colon (J,<)\models s<t_{\alpha}\}\big)=\k$.
    
    Let $$ I_{\al}^a=\{ t\in I_{\al}\colon t<t_{\al}\},$$
    $$I^b_{\al}=\{ t\in I_{\al +1} \colon t\in I_{\al}^a \vee tEt_{\al} \}, $$ $$I_{\al}^c=\{ t\in I_{\al} \colon t> t_{\al}\}.$$
     
    Note that $I_{\al}=I^a_{\al}+I^c_{\al}$ and
    $I_{\al+1}=I^b_{\al}+I^c_{\al}$. Recalling Fact (C), it is now enough to
    show that
    $I^a_{\al}\subnice I^b_{\al}$. Identifying isomorphic orders and using $(\a)$,
    one has that $\unl\k$ is unbounded in $I_{\al}^a$ and $I_{\al}^b=I^a_{\al}+
    (I_{\al}^b\setminus I_{\al}^a)$ so by Fact (B), $I_{\al}^a\subnice I_{\al}^b$
    as required.
\end{PROOF}

Of the five facts, we prove (A), (B), and (E) as (C) and (D) are obvious.

\begin{proof}[{Proof of Fact (A)}]
    Since $\k$ is measurable, there is a $\k$-complete
    uniform ultra-filter $D$ on $\k$ (see \cite{Je03}). For every linear order $J_0$
    (or $J^*_0$) there is
    $\Op_{I,D}(-)$, the iteration of $I$ ultra-powers $(-)^{\k}/D$,
    ordered in the order $J_0$
    (or $J^*_0$), giving the required embedding (use Observation \ref{w17}).
\end{proof}

\begin{proof}[{Proof of Fact (B)}]
    Since $\underline{\k}\subseteq I$ and using
    Fact (A), we know that

    \begin{enumerate}
        \item[$\bullet_{1}$] let $d_{0}$ be the identity map from $\underline{\kappa}$ into $I$,

        \item[$\bullet_{2}$] let $d_{1}$ be the canonical map from $\underline{\kappa}$ into $\Op(\underline{\kappa})$, which exists by properties of $\Op,$

        \item[$\bullet_{3}$] let $d_{2}$ be the embedding of $J_{0}$ into $\Op(\underline{\kappa})$ as in the choice of $\Op$, so $\rang(d_{1})$ is below $\rang(d_{2}),$

        \item[$\bullet_{4}$] let $d_{3}$  be the canonical embedding from $\Op(\underline{\kappa})$ into $\Op(I)$ by 
        lifting 
        (really is the identity), 

        \item[$\bullet_{5}$] let $d_{4}$ be the canonical embedding of $I$ into $\Op(I)$, which exists by properties of $\Op$ extending $d_{1}$ by the properties of $\Op,$

        \item[$\bullet_{6}$] So $\rang(d_{4} \rest \kappa)$ is unbounded in $\rang(d_{4})$ in the order $\Op(I)$,

        \item[$\bullet_{7}$] $\rang(d_{4}) \subseteq \Op(I)$ is below $d_{3} \circ d_{2}'' \in J_{0}.$ 
    \end{enumerate}
    
    
    
    Chasing through the diagram, we obtain the required embedding.     So we are done. 
\end{proof}

\begin{proof}[{Proof of Fact (E)}]
     Apply Observation \ref{X1.11} to the chain
    $\lng I_\alpha\colon\al\le \del\rng$. 
\end{proof}

So we are done proving Criterion \ref{X2.6}. 

\begin{Fact}\label{X2.7}
    Suppose that $\lambda \ge \k$.  There exist a linear order $(I, <_I)$ of power $\lambda$ and a sequence $\langle A_i \subseteq I \colon i \le \lambda \rangle$ of pairwise disjoint subsets of $I$, each of power $\k$
    such that $I = \underset{i \le \lambda} \cup  A_i$ and,
    
    $(\ast)$ if $\lambda \in X \subseteq \lambda \dotplus 1$, then $I \uhr \underset{i \in X} \cup A_i \underset{\text{nice}}  \subseteq I$.
\end{Fact}

\begin{PROOF}{\ref{X2.7}}
    Let $I = (\lambda \dotplus 1) \times \k$ and define $<_I$ on $I\colon (i_1, \alpha_1) <_I (i_2, \alpha_2)$ iff $i_1 < i_2$ or $(i_1 = i_2$ and $\alpha_1 > \alpha_2)$.  For each $i \le \lambda$, let $A_i = \{ i \} \times \k$. Let us check $(\ast)$ of Criterion \ref{X2.6}: suppose that $\lambda \in X \subseteq \lambda + 1$.  Write $I_X = I
    \uhr \left( \bigcup_{i \in X} A_{i} \right)$. 
    To show that $I_X \underset{\text{nice}}  \subseteq I$,  we can assume without loss of generality that $I_{X} \neq I$ and then one employs  Criterion \ref{X2.6}. 
    Consider $t \in I - I_X$, say $t = (i, \alpha)$ (note that $\alpha < \k$ and $i < \lambda$, since $\lambda \in X$) and $i \notin X$.  Let $j = \min (X -i)$; note that $j$ is
    well-defined, since $\lambda \in X -i$, and $j \neq i$. For every $\beta < \k$, one has $t <_I (j, \beta)$ and $(j, \beta) \in I_X$.
    Also if $s \in I_X$ and $t <_I s$, then for some $\beta < \k  (j, \beta) <_I s$.  Thus $\langle (j, \beta)\colon  \beta < \k \rangle$ is a
    cofinal sequence in $( I_X \uhr \{ s \in I\colon  t <_I s \})^\ast$.  By the criterion, $I_X \underset{\text{nice}} \subseteq I$.
\end{PROOF}

\begin{Theorem}\label{X2.8}
    Suppose that $\k = \cf(\delta) \le \delta < \lambda$. Then $\EM(\delta) \underset{\text{nice}}  \leq \EM(\lambda)$.
\end{Theorem}

\begin{PROOF}{\ref{X2.8}}
    By Fact (B) of Criterion \ref{X2.6}, one has that $\d\subnice \lambda$; so by Theorem \ref{X2.5},
    $\EM(\d) {\underset{\text{nice}}  \leq} \EM(\lambda)$.
\end{PROOF}

Now let us turn to the main theorem of this section.

\begin{Theorem}\label{X2.9}
    Suppose that $\gk$ is categorical in the regular
    cardinal $\lambda > \chi$.
    \underline{Then} ${\gk}_{ < \lambda}$ has the amalgamation property.
\end{Theorem}

\begin{PROOF}{\ref{X2.9}}
    Suppose that $\gk_{< \lambda}$ fails AP.  Note that $\| \EM(\lambda) \| = \lambda$.  Apply Theorem \ref{X2.4} to find $M^\ast \in K_\lambda$ and $\langle M_i \colon i < \lambda \rangle$ satisfying Theorem
    \ref{X2.4}(1) and Theorem \ref{X2.4} (2).  Since $\gk$ is
    $\lambda$-categorical, $M^\ast \cong E M (\lambda)$, so without loss of generality $\EM(\lambda)
    = \underset{i < \lambda} \cup M_i$.  $C = \{ i < \lambda\colon  M_i = \EM(i) \}$ is a club of $\lambda$.  Choose $\delta \in C, \cf(\delta) = \k$.  By Theorem \ref{X2.8}, $\EM(\delta) \underset{\text{nice}}  \leq \EM(\lambda)$, so $M_\delta \underset{\text{nice}}  \leq M^\ast$.  But of course by Theorem \ref{X2.4}(2) $M_\delta \underset{\text{nice}}  
    \nleq 
    M^\ast$-contradiction.
\end{PROOF}

The last theorem of this section applies to singular cardinals as well. 

\begin{Theorem}\label{X2.10}
    Suppose that $K$ is categorical in
    $\lambda > \chi$ (notice that $\lambda$ is not necessarily regular).  Then:
    
    (1) $K$ has a model $M$ of power $\lambda$ such that if  $N \leq_{\gk} M$ and
    $\| N \| < \lambda$, then there
    exists $N'$ such that:
    
    \begin{enumerate}
        \item[$(\alpha)$] $N \leq_{\gk} N' \leq_{\gk} M,$
        
        \item[$(\beta)$] $\| N' \| = \| N \| + \chi,$
        
        \item[$(\gamma)$] $N' \underset{\text{nice}}  \leq M.$
    \end{enumerate}
    
    (2) $K$ has a model $M$ of power $\lambda$ and an expansion $M^+$ of $M$ by at most $\chi$
    functions such that if $N^+ \subseteq M^+$, then $N^+ \uhr \tau \underset{\text{nice}}  \leq M$.
\end{Theorem}

\begin{PROOF}{\ref{X2.10}} 
    (1) Let $\langle I, \langle A_i \colon  i \le \lambda \rangle \rangle$
    be as in Fact \ref{X2.7}. 
    Let $M = \EM(I)$.  Suppose that $N \leq_{\gk} M,  \| N \| < \lambda$.
    Then there exists $J \subseteq I,
    |J| < \lambda$ such that $N \subseteq \EM(J)$ so by Fact \ref{X2.7} there exists $X \subseteq \lambda + 1$ such that $\lambda \in X,$ $\vert X \vert < \lambda$ and $J \subseteq \underset{i \in X} \cup A_i$. 
    
    Note that $\left| \underset{i \in X} \cup A_i \right| \le \vert X \vert \cdot \kappa < \lambda$.
    Now $N' = \EM(I \uhr \underset{i \in X} \cup A_i)$
    is as required, since
    $I \uhr \underset{ i \in X} \cup A_i \underset{\text{nice}}  \leq I$ and so by Theorem \ref{X2.5} 
    $\EM(I \uhr \left( \bigcup_{i \in X} A_{i} \right)) \underset{\text{nice}}  \leq \EM(I)$.  This proves (1).
    
    (2) We expand $M = \EM(I)$ with skeleton $\langle a_{s} \colon s \in I \rangle$ as follows:

    \begin{enumerate}
        \item[(a)] by all functions of $\EM'(I),$
        
        \item[(b)] by the unary functions $f_\ell(\ell<n)$ which are chosen as
        follows: we know that for each $b \in M$ there is $\sigma_b$ an $\tau_1$-term ($\tau_1$-the vocabulary of $\EM'(I)$) and $t(b,0) < t(b,1) < \ldots < t(b,n_{\sigma_b}-1)$ from $I$ such that $$b=\sigma_b(a_{t(b,0)}, a_{t(b,1)}, \dots, a_{t(b, n_{\sigma_{b}}-1)})$$
        
        (it is not unique, but we can choose one; really if we choose it with $n_b$ minimal it is almost unique). We let $$ f_\ell(b) =
        \begin{cases} a_{t((b,\ell))}, &\text{ if } \ell<n_{\sigma_b},\\
        b, &\text{ if } \ell\ge n_{\sigma_b}.\\ \end{cases}
        $$
        
        \item[(c)] by unary functions $g_\al$, $g^\al$ for $\al<\k$ such that
        \underbar{if} $t  < s$ are in $I$, $\al = \otp[(t,s)^*_I]$ then
        $g^\al(a_t)=a_s$, $\bigvee\limits_{\beta< \kappa}g_\beta(a_s)=a_t$
        (more formally $g^\alpha(a_{(i, \beta)})=a_{(i, \beta+\alpha)}$ and
        $g_\alpha(a_{(i, \beta)}) = a_{(i, \alpha)}$)
        in the other cases $g^\al(b)=b$, $g_\al(b)=b$.
        
        \item[(d)] by individual constants $c_\al=a_{(\lambda,\al)}$ for $\al<\k$.
    \end{enumerate}
    
    Call the expanded model $M^{+}$. Now suppose $N^+$ is a sub-model of $M^+$ and $N$ its $\tau$-reduct.
    Let $J\eqdf\{t\in I\colon a_t\in N\}$, now $J$ is a subset of $I$ of cardinality $\le||N||$ as for $t\not=s$ from $J$, $a_t\not= a_s$. Also if $b\in N$ by clause (b), $a_{t(b,\ell)}\in N$ hence $b\in \EM(J)$; on the other hand if $b\in \EM(J)$ then by clause (b) we have $b\in N$; so we can conclude $N=\EM(J)$. So far this holds for any linear suborder of $I$.
    
    By clause (c) $J=\bigcup\limits_{i\in X}A_i$ for some $X\subseteq \lam+1$, and by clause (d), $\lambda \in X$.
    
    Now $\EM(J) \underset{\rm{nice}}{\leq} \EM(I) = M$ by Fact \ref{X2.7}.
\end{PROOF}

\newpage 

\section{Towards removing the assumption of regularity from the existence of universal extensions}\label{3}

In \S 2 we showed that $\gk_{< \lambda}$ has the amalgamation property when $\gk$ is categorical in the regular cardinal $\lambda > \chi$.  We now study the situation in which $\lambda$ is not assumed to be regular.

Our problem is that while we know that most sub-models of $N\in K_\lam$ sit well in $N$ (see Theorem \ref{X2.10}(2)) and that there are quite many $N\in K_{<\lam}$ which are amalgamation bases, our difficulty is to get those things together:  constructing $N\in K_{\lam}$ as $\bigcup\limits_{i<\lam}N_i$, $N_i\in K_{<\lam}$ means $N$ has $\leq_{\gk}$-sub-models not included in any $N_i$.

Recall we are assuming Hypothesis \ref{b2}. 

\begin{Theorem}\label{X3.1} 
    Suppose that $\gk$ is categorical in $\lambda$ and $\chi \le \theta < \lambda$.  If $\langle M_i \in K_\theta \colon i < \theta^+ \rangle$ is an increasing continuous $\leq_{\gk}$-chain, then:
    $$ \left \{ i < \theta^+  \colon M_i \underset{\text{nice}}  \leq (\cup_{j < \theta^{+}} M_{j})  \right \} \in D_{\theta^+}. $$
\end{Theorem}

\begin{Remark}\label{X3.1A}\ 
    
    (1) We cannot use Theorem \ref{X2.10}(1) as possibly $\lambda$ has cofinality $< \chi$.
    
    (2) Recall that $D_{\theta^+}$ is the closed unbounded filter on $\theta^+$.
\end{Remark}

\begin{PROOF}{\ref{X3.1}}
    Write $M_{\theta^+} = \underset{i < \theta^+} \cup M_i$. Choose an operation $\Op$ such that for all $i < \theta^+, \, \| \Op(M_i) \| \ge \lambda$.
    Let $M^\ast_i = \Op(M_i)$, hence $M_{i} \underset{\mathrm{nice}}{\leq} M_{i}^{\ast}.$ Applying Fact \ref{X1.2}  for non-limit ordinals, Fact \ref{X1.1} for limit ordinals, one finds inductively an increasing continuous $\leq_{\gk}$-chain $\langle N_i \colon  i \le \theta^+ \rangle$ such that for $i < \theta^+, M_i \leq_{\gk} N_i \leq_{\gk} M^\ast_i, \| N_i \| = \lambda,$ so $M_{i} \underset{\mathrm{nice}}{\leq} N_{i}$ and $
    N_{\theta^+} = \underset{i < \theta^+} \cup N_i$.  Note that $\| N_{\theta^+} \| = \theta^+ \cdot \lambda = \lambda.$ 
    
    Since $\gk$ is $\lambda$-categorical, $N_{\theta^+} \cong \EM(I)$ where Fact \ref{X2.7}  furnishes $I$ of power $\lambda$.  By Theorem \ref{X2.10}(2), there is an expansion $N^+_{\theta^+}$ of $N_{\theta^+}$ by at most $\k + | \tau_{\gk} |$ functions such that if $A \subseteq |N^+_{\theta^+}|$ is closed under the functions of $N^+_{\theta^+}$, then $(N^+_{\theta^+} \uhr 
      \tau _{\gk}
    ) \uhr A \underset{\text{nice}}  \leq N_{\theta^+}$.
    
    Choose a set $A_i$ and an ordinal $j_i$, by induction on $i < \theta^+$, satisfying:
    
    (1) $A_i \subseteq |N_{\theta^+}|,  |A_i| \le \theta;\langle A_i\colon  i < \theta^+ \rangle$ is continuous increasing, 
    
    (2) $\langle j_i \colon  i < \theta^+ \rangle$ is continuous increasing, 
    
    (3) $A_i$ is closed under the functions of $N^+_{\theta^+}$, 
    
    (4) $A_i \subseteq |N_{j_{i + 1}} |$,
    
    (5) $|M_i| \subseteq A_{i+1}$.

    This is possible: for zero, let $A_{0} \coloneqq \emptyset,$ $j_{0} \coloneqq 0$ and for limit ordinals unions work; for $i +1$ choose $j_{i+1}$ to satisfy (2) and (4), and $A_{i+1}$ to satisfy (1), (3) and (5).
    
    By (2), $C = \{ i < \theta^+ \colon  i$ is a limit ordinal and $j_i = i \}$ is a club of $\theta^+$ i.e. $C \in D_{\theta^+}$.
    
    Fix $i \in C$.  Note that $|M_i | \subseteq A_i$ and $A_i \subseteq |N_i|$ (since $|M_i| = \underset{j < i}  \cup | M_j | \subseteq \underset {j < i} \cup
    A_{j +1} = A_i = \underset{i' < i} \cup A_{i'} \subseteq \underset{i' < i} \cup | N_{ j_{i' + 1}} | = N_{j_i} = N_i$ (using (5), (1), (4), (2) and
    $j_i = i$)) and so $M_i \leq_{\gk} (N^+_{\theta^+} \uhr 
    \tau _{\gk} 
    )\uhr A_i \leq_{\gk} N_i \leq_{\gk} M^\ast_i = \Op(M_i)$, so that $M_i \underset{\text{nice}} \leq (N^+_{\theta^+} \uhr 
     \tau _{\gk}
    ) \uhr A_i$. However by (3) and the choice of $N_{\theta^+}$ and $N^+_{\theta^+}$ one has also that $(N^+_{\theta^+} \uhr 
    \tau _{\gk}
    ) \uhr A_i \underset{\text{nice}}  \leq
    N_{\theta^+}$.  So by transitivity of $\underset{\text{nice}}  \leq$, one obtains
    $M_i \underset{\text{nice}}  \leq N_{\theta^+}$.
    
    Finally remark that $M_{\theta^+} \leq_{\gk} N_{\theta^+}$ since
    $M_i \underset{\text{nice}} \leq N_i \leq_{\gk} N_{\theta^+}$ for
    every $i < \theta^+$. Hence $i \in C \Rightarrow M_{i} \leq_{\gk} M_{\theta^{+}} \leq_{\gk} N_{\theta^{+}},$ so recalling $i \in C \Rightarrow M_{i} \underset{\mathrm{nice}}{\leq} N_{\theta}$ we have $i \in C \Rightarrow M_{i} \underset{\mathrm{nice}}{\leq} M_{\theta^{+}};$ so $C \subseteq \left \{ i < \theta^+\colon  M_i \underset{\text{nice}} \leq
    M_{\theta^+} \right \} \in D_{\theta^+}$.
\end{PROOF}

\begin{Definition}\label{X3.2}
    Suppose that $\theta \in [\chi, \lambda )$ and $M \in K_{\theta}$. $M$ is nice iff whenever $M \leq_{\gk} N \in K_\theta$, then $M \underset{\text{nice}}  \leq N$.  (The analogous $\leq_{\gk}$-elementary embedding definition runs: $M$ is nice iff whenever $f \colon M \underset{\gk}  \rightarrow N \in K_\theta$ then $f \colon M \underset{\text{nice}}  \rightarrow N$).
\end{Definition}

\begin{Theorem}\label{X3.3}
    Suppose that $\gk$ is categorical in $\lambda$ and $M \in K_\theta, \theta \in [\chi, \lambda)$.  Then there exists $N \in K_\theta$ such that $M \leq_{\gk} N$ and $N$ is nice.
\end{Theorem}

\begin{PROOF}{\ref{X3.3}} 
    Suppose otherwise.  We'll define a continuous increasing $\leq_{\gk}$-chain
    $\langle M_i \in K_\theta\colon  i < \theta^+ \rangle$ such that for $j < \theta^+$:
    
    \begin{enumerate}
        \item[$(*)_{j}$] $M_j \underset{\text{nice}} \nleq M_{j +1}$.
    \end{enumerate}
    
    For $i = 0$, put $M_0 = M$;  if $i$ is a limit ordinal, put $M_i = \underset{j < i}
    \cup M_j$;  if $i = j +1$, then since Theorem \ref{X3.3} is assumed to fail, $M_{j +1}$
    exists as required in $(\ast)_j$ (otherwise $M_j$ works as $N$ in Theorem \ref{X3.3}).  But now
    $\langle M_i \colon  i < \theta^+ \rangle$ yields a contradiction to Theorem \ref{X3.1}, since $C = \{ i < \theta^+ \colon  M_i  \underset{\mathrm{nice}}{\leq} \underset{ j < \theta^+} \cup M_j \} \in D_{\theta^+}$ by Theorem \ref{X3.1} so that choosing $j$ from $C$ one has $M_j \underset{\text{nice}} \leq M_{j +1}$ by Observation \ref{X1.10}, contradicting $(\ast)_j$.
\end{PROOF}

\begin{Theorem}\label{X3.4}
    Suppose that $\gk$ is categorical in $\lambda$ and $\theta \in [\chi, \lambda)$.  If $M \in K_\theta$ is nice and $f\colon  M \underset{\gk}  \rightarrow N \in K_{\le \lambda}$, then $f\colon  M \underset{\text{nice}}  \rightarrow N$.
\end{Theorem}

\begin{PROOF}{\ref{X3.4}} 
    Choosing an appropriate $\Op$ and using Fact \ref{X1.2} one finds $N_1$ such that $N \leq_{\gk} N_1$ and $\| N_1 \| = \lambda$. Find $M'_1 \underset{\text{nice}}  \leq N_1$ by Theorem  \ref{X2.10}(2) such that $\mathrm{rng}(f) \subseteq |M'_1|,  \| M'_1 \| = \theta$.  So $M'_1 \leq_{\gk} N_1$ and therefore $ N_{1} \, {\rest} \,  \mathrm{rng}(f) \leq_{\gk} M'_1$. Recall  $M$ is nice, so $f \colon  M \underset{\text{nice}}  \rightarrow M'_1$. Now $M'_1 \underset{\text{nice}}  \leq N_1$, therefore $f\colon  M \underset{\text{nice}}  \rightarrow N_1$. So there are $\Op$ and $g \colon N_1 \underset{\gk} \rightarrow \Op(M)$ satisfying $gf = f_{\Op}$. Since $N \leq_{\gk} N_1$ it follows that $f \colon  M \underset{\text{nice}} \rightarrow N$ as required.
\end{PROOF}

\begin{Corollary}\label{X3.5}
    Suppose that $M \in K_\theta$ is nice, $\theta \in [\chi, \lambda)$. Then $M$ is an a.b. in $\gk_{\le \lambda}$ i.e. if $f_i \colon M \underset{\gk} \rightarrow M_i, M_i \in K_{\le \lambda} (i = 1,2)$, then there exists an amalgam $N \in K_{\le \lambda}$ of $M_1, M_2$ over $M$ w.r.t. $f_1, f_2$.
\end{Corollary}
    
\begin{PROOF}{\ref{X3.5}}
    By Definition \ref{X3.4} $f_i \colon  M \underset{\text{nice}}  \rightarrow
    M_i  (i =1, 2)$.  Hence by Lemma \ref{X2.1} 
    there is an amalgam $N \in K_{\le \lambda}$ of $M_1, M_2$ over $M$
    w.r.t. $f_1, f_2$.
\end{PROOF}

\begin{Definition}\label{X3.6}
    Suppose that $\theta \in [\chi, \lambda)$ and  $\partial$ is a cardinal.
    
    (1) A model $M \in K_\theta$ is $\partial$-universal iff for every
    $N \in K_\partial$, there exists an $\leq_{\gk}$-elementary embedding $f\colon  N \underset{\gk} \rightarrow M$. We say $M$ is universal iff $M$ is $\| M \|$-universal.
    
    (2) A model $M_2 \in K_{\theta}$ is $\partial$-universal over the model $M_1$ (and one writes $M_1 \underset{\partial-\rm{univ}} \preceq M_2$) \underline{iff} $M_1 \leq_{\gk} M_2$ and whenever $M_1 \leq_{\gk} M'_2 \in K_{\partial}$, then there exists an $\leq_{\gk}$-elementary embedding $f\colon  M'_2 \underset{\gk}  \rightarrow M_2$ such that $f \uhr M_1$ is the identity.  (The embedding version runs:  there exists $h\colon  M_1 \underset{\gk} \rightarrow M_2$ and whenever $g \colon  M_1 \underset{\gk} \rightarrow M'_2 \in K_{ \partial }$, then there exists $f\colon  M'_2 \underset{\gk} \rightarrow M_2$ such that $f g = h$.) $M_2$ is universal over $M_1 \ (M_1 \underset{\text{univ}} \preceq M_2)$
    iff $M_2$ is $\|M_2 \|$-universal over $M_1$.
    
    (3) A model $M_2$ is $\partial$-universal over $M_1$ in $M$ \underline{iff} $M_1 \leq_{\gk} M_2 \leq_{\gk} M$, $||M_1||\le\partial$ and whenever $M'_2 \in K_{\partial}$ and $M_1 \leq_{\gk} M'_2 \leq_{\gk}M$, then there exists an $\leq_{\gk}$-elementary embedding $f\colon  M'_2 \underset{\gk} \rightarrow M_2$ such that $f \uhr M_1$ is the identity.  $M_2$ is universal over $M_1$ in $M$ iff $M_2$ is $\| M_2 \|$-universal over $M_1$ in $M$.
    
    (4) $M_2$ is weakly $\partial$-universal over $M_1$ (written $M_1 \underset{\partial- \rm{wu}} \prec M_2$) \underline{iff} $M_1 \leq_{\gk} M_2\in K_\partial$ and whenever $M_2 \leq_{\gk} M'_2 \in K_{\partial}$, then there exists an $\leq_{\gk}$-elementary embedding $f\colon  M'_2 \underset{\gk}  \rightarrow M_2$ such that $f \uhr M_1$ is the identity.  (The embedding version is:  there exists $h\colon  M_1 \underset{\gk}  \rightarrow M_2$ and whenever $g \colon  M_2 \underset{\gk} \rightarrow M'_2 \in K_{\partial}$, then there exists $f \colon  M'_2 \underset{\gk} \rightarrow M_2$ such that $h = f g h$ (written $h\colon  M_1 \underset{\partial-\rm{wu}} \rightarrow M_2$)). We say  $M_2$ is weakly universal over $M_1 (M_1 \underset{\rm{wu}} \preceq M_2)$ iff $M_2$ is $\|M_2 \|$-weakly universal over $M_1$.
\end{Definition}


\begin{Remark}\label{X3.6A}
     In $\gk_{<\lam}$, if $M_1$ is an a.b., then weak universality over
    $M_1$ is equivalent to universality over
    $M_1$.
\end{Remark}

\begin{PROOF}{\ref{X3.6A}}
    Suppose that $h\colon  M_1 \underset{\rm{wu}} \rightarrow M_2$ and $g\colon  M_1 \underset{\gk} \rightarrow M'_2 \in K_{\| M_2 \| }$. Since $M_1$ is an a.b. there exist a model $N$ and $h'\colon  M_2 \underset{\gk} \rightarrow N, \; g'\colon  M'_2 \underset{\gk} \rightarrow N$ satisfying $h' h = g'g$.  By Fact \ref{X1.2} without loss of generality $\| N \| = \| M_2 \|$.  Since $M_2$ is weakly universal over $M_1$, there exists $h^{''}\colon  N \underset{\gk} \rightarrow M_2, \; h = h^{''} h' h$. Let $f = h^{''} g' \colon  M'_2 \to   M_2$, and note that $fg \uhr M_1 = h^{''} g'g = h^{''} h' h = h$, so that $M_2$ is universal
    over $M_1$.
\end{PROOF}

\begin{Remark}\label{X3.6B}
    For any model $M$, universality over $M$ implies weak universality over $M$.
\end{Remark}

\begin{Lemma}\label{X3.7} 
    Suppose that $\gk$ is categorical in $\lambda, \theta \in [\chi, \lambda)$. If $M \in K_\theta$ and $M \leq_{\gk} N \in K_\lambda$, then there exists $M^+ \in K_\theta$ such that:
    
    \begin{enumerate}
        \item[(a)] $M \leq_{\gk} M^+ \leq_{\gk} N,$
        
        \item[(b)] $M^+$ is universal over  $M$ in $N.$
    \end{enumerate}
\end{Lemma}

\begin{PROOF}{\ref{X3.7}}
    We choose $I$ such that:
    
    \begin{enumerate}
        \item[$(\ast)$] 
        
        \begin{enumerate}
            \item[(a)]  $I$ is a linear order of cardinality $\lambda,$
            
            \item[(b)] if $\partial \in [\aleph_0,\lam)$, $J_0\subseteq I$, $|J_0|= \partial$ then there is $J_1$ satisfying $J_0\subseteq J_1\subseteq I$, $|J_1|= \partial$, and for every $J^*\subseteq I$ of cardinality $\le \partial$ there is an order-preserving (one to one) mapping from $J_0\cup J^*$ into $J_0\cup J_1$ which is the identity on $J_0. $
        \end{enumerate}
    \end{enumerate}

    Essentially the construction follows Laver \cite{Lv71}  and  \cite[Appendix]{Sh:220}, see more in \cite{Sh:E62}; but for our present purpose let $I=({\;}^{\om>}\lam,<_{\ell \rm{e x}})$; given $\theta$ and $J_0$ we can increase
    $J_0$ so without loss of generality $J_0={\;}^{\om>}A$, $A\subseteq \lam$, $|A|=\theta$. Define an equivalence relation $E$ on $I\setminus J_0$: $\eta E \nu\Leftrightarrow (\forall \rho \in J_0)(\rho<_{\ell \rm{ex}}\eta \equiv \rho<_{\ell \rm{ex}}\nu)$, easily it has $\le\theta$ equivalence classes, so let $\{\eta^*_i \colon i<i^*\le\theta\}$
    be a set of representatives each of minimal length, so
    $\eta_{i}^*\uhr(\lg\eta^*_i-1) \in J_0$, $\eta^*_i(\lg \eta^*_i-1)\in \lam\setminus A$.
    
    Let $J_1=I\cup\{\eta^*_i{\;}^{\hat{\; }}\nu\colon  \nu\in{\;}^{\om>}\theta$ and $i<i^*\}$, so clearly
    $J_0\subseteq J_1\subseteq I$, $|J_1 \vert =\theta$. Suppose
    $J_0\subseteq J\subseteq I$, $|J|\le\theta$, and we should find the required embedding $h$. As before without loss of generality $J={\;}^{\om>}B$, $|B|=\theta$ and $A\subseteq B$. Now $h\uhr J_0=\rm{id}_{J_0}$ so it is enough to define $h\uhr
    (J_1\cap(\eta^*_i/E))$, hence it is enough to embed $J_1\cap(\eta^*_i/E)$ into $\{\eta_1^*{\;}^{\hat {\; }} \nu\colon  \nu\in{\;}^{\om>}\theta\}$ (under $<_{\ell \rm{ex}}$).
    
    Let $\gam = \otp(B)$, so it is enough to show
    $({\;}^{<\om}\gam,<_{\ell \rm{ex}})$ can be embedded into ${\;}^{\om>}\theta$, where of course $|\gam|\le\theta$. This is proved by induction on $\gam$.
    
    Since $\gk$ is $\lambda$-categorical and $\EM(I)$ is a model of $\gk$ of power $\lambda$, there is an isomorphism $g$ from $\EM(I)$ onto $N$. It follows from $(\ast)$ that $M^+ = g^{''}(\EM(J)) \in K_\theta$ satisfies (1) and (2). (Analogues of (1) and (2) are checked in more detail in the course of the proof of Corollary \ref{X3.11}.)
\end{PROOF}

\begin{Lemma}\label{X3.8} 
    Suppose that $\gk$ is categorical in
    $\lambda, \theta \in [\chi, \lambda)$, and $\langle M_i \in K_\theta \colon i < \theta^+ \rangle, \, \langle N_i \in K_\lambda \colon i < \theta^+\rangle$ are continuous $\leq_{\gk}$-chains such that for every $i < \theta^+$ we have $M_i \leq_{\gk}  N_i$.  Then there exists $i (\ast) < \theta^+$ such that $(i (\ast),\; \theta^+) \subseteq C \coloneqq \{ i < \theta^+\colon  M_{i +1}$ can be $\leq_{\gk}$-elementarily embedded into $N_i$ over $M_0 \}$.
\end{Lemma}

\begin{PROOF}{\ref{X3.8}}
    Apply Lemma \ref{X3.7} for $M_0 \in K_\theta$ and $N_{\theta^+} = \underset{i < \theta^+} \cup N_i \in K_\lambda$ (noting that $M_0 \leq_{\gk}  N_0 \leq_{\gk} N_{\theta^+}$) to find $M^+ \in K_\theta$ such that  $M_0 \leq_{\gk} M^+ \leq_{\gk} N_{\theta^+}$ and $M^+$ is universal over $M_0$ in
    $N_{\theta^+}$.

    For some $i (\ast) < \theta^+, M^+ \subseteq N_{i (\ast)}$ and so $M^+ \leq_{\gk} N_{i (\ast)}$.  If $i \in (i (\ast),  \theta^+)$, then $M_{i +1} \in K_\theta$ and $M_0 \leq_{\gk} M_{i +1} \leq_{\gk} N_{i +1} \leq_{\gk} N_{\theta^+}$, so there is an $\leq_{\gk}$-elementary embedding $f \colon  M_{i +1} \underset{\gk} \rightarrow M^+$ and $f \uhr M_0$ is the identity.  Now $M^+ \leq_{\gk} N_{i (\ast)} \leq_{\gk} N_i$, so
    $f\colon  M_{i +1} \underset{\gk} \rightarrow N_i$. Hence $(i (\ast), \; \theta^+) \subseteq C$ as required.
\end{PROOF}

\begin{Theorem}\label{X3.9}
    Suppose that $\gk$ is categorical in $\lambda, \theta \in [\chi, \lambda),  M \in K_\theta$.  \underline{Then}
    there exists $M^+ \in K_\theta$ such that:
    
    \begin{enumerate}
        \item[$(\aleph)$]  $M \leq_{\gk} M^+$  and  $M^+$ is nice,
        
        \item[$(\beth)$] $M^+$ is weakly universal over $M.$
    \end{enumerate}
\end{Theorem}

\begin{PROOF}{\ref{X3.9}}
    Define by induction on $i < \theta^+$ continuous $\prec_{\gk}$-chains
    $\langle M_i \in K_\theta \colon  i < \theta^+ \rangle, \, \langle N_i \in
    K_\lambda\colon  i < \theta^+ \rangle$ such that:
    
    (0) $M_0 = M,$
    
    (1) $M_i \leq_{\gk} N_i,$
    
    (2) if $(\ast)_i$ holds, then $M_{i +1}$ cannot be $\leq_{\gk}$-elementarily embedded into $N_i$ over $M_0$,
    where $(\ast)_i$ is the statement:
    
    \begin{enumerate}
        \item[$(*)_{i}$] there are $M' \in K_\theta$ and $N' \in K_\lambda$ such that $M_i \leq_{\gk} M'$, $N_i \leq_{\gk} N'$, $M' \leq_{\gk} N'$ and $M'$ cannot be $\leq_{\gk}$-elementarily embedded into $N_i$ over $M_0,$
    \end{enumerate}
    
    (3) $M_{i +1} \underset{\text{nice}}  \leq N_{i +1}$.
    
    This is possible.  $N_0$ is obtained by an application of Fact \ref{X1.2} to an appropriate $\Op(M_0)$ of power at least $\lambda$.  At limit stages, continuity dictates that one take unions.  Suppose that $M_i, N_{i}$ have been defined. If $(\ast)_i$ does not hold, by Theorem \ref{X2.10}(2) there is $M''\in K_{\theta}$, $M_i \leq_{\gk} M''
      \underset{\rm nice}  \leq 
    N_i$. Let $M_{i +1} = M'', \; N_{i +1} = N_i$. If $(\ast)_i$ does hold for $M', N'$, let $N_{i +1} = N'$; note that by Theorem \ref{X2.10}(2) there exists $M^{''} \in K_\theta ,  M' \leq_{\gk} M^{''} \underset{\text{nice}}  \leq N'$; now let $M_{i +1} = M^{''}$.  Note that in each case, (3) is satisfied.
    
    Find $i (\ast) < \theta^+$ and $C$ as in Lemma \ref{X3.8} and choose $i \in C$. By (1), we have  $M_{i +1} \leq_{\gk} N_{i +1}$ so by Lemma  \ref{X3.7} there exists $M^- \in K_\theta$ such that $M_{i +1} \leq_{\gk} M^- \leq_{\gk} N_{i +1}$ and $M^-$ is weakly universal over $M_{i +1}$ in $N_{i +1}$.  By Theorem \ref{X3.3} one can find $M^+ \in K_\theta$ such that
    $M^- \leq_{\gk} M^+$ and $M^+$ is nice.  So $M^+$ satisfies $(\aleph)$. It remains to show that $M^+$ is weakly universal over $M$.  Suppose not and let $g\colon  M^+ \underset{\gk}  \rightarrow M^\ast \in K_\theta$ where $M^\ast$ cannot be $\leq_{\gk}$-elementarily embedded in $M^+$ over $M$ hence cannot be $\leq_{\gk}$-elementarily embedable in $M^-$ over $M$, hence in $N_{i+1}$ over $M$. $M_{i +1} \leq_{\gk} M^\ast\in K_\theta$ and by (3) $M_{i +1} \underset{\text{nice}}  \leq N_{i +1}\in K_\lambda$, so by
    2.1 there is an amalgam $N^\ast \in K_\lambda$ of $M^\ast, \; N_{i +1}$. The existence of $M^\ast, N^\ast$ implies that $(\ast)_{i +1}$ holds since $M^\ast$ cannot be $\leq_{\gk}$-elementarily embedded into $N_{i +1}$ over $M_0$, hence $M_{i +2}$ cannot be $\leq_{\gk}$-elementarily embedded into $N_{i +1}$ in contradiction to
    the choice of $i$ as by Lemma \ref{X3.7} $i +1$ is in $C$.
\end{PROOF}

\begin{Corollary}\label{X3.10} 
    If $\gk$ is categorical in $\lambda, \theta \in [\chi, \lambda)$ and $M \in K_\theta$ is an a.b. (e.g. $M$ is nice, see \ref{X2.1}), then there exists $M^+ \in K_\theta$ such that:
    
    \begin{enumerate}
        \item[$(\aleph)$]   $ M \leq_{\gk} M^+$ and $M^+$ is nice,
        
        \item[$(\beth)$] $M^+$ is universal over $M$.
    \end{enumerate}
\end{Corollary}

\begin{PROOF}{\ref{X3.10}} 
    By Theorem \ref{X3.9} and Remark \ref{X3.6A}.
\end{PROOF} 

\begin{Corollary}\label{X3.11}
    Suppose that $\gk$ is categorical in $\lambda$ and $\theta\in [\chi, \lambda)$. \underline{Then} there is a nice universal model $M \in K_\theta$.
\end{Corollary}
    
\begin{PROOF}{\ref{X3.11}}
    By \ref{X3.3} it suffices to find a universal model of power $\theta$, noting
    that universality is preserved
    under $\leq_{\gk}$-elementary extensions in the same power. As in the proof of \ref{X3.7}, there is a linear order $(I, <_I)$ of power $\lambda$ and $J \subseteq I, |J| = \theta$, such that: 
    
    \begin{enumerate}
        \item[$(*)$] $ ( \forall J' \subseteq I)$ (if $| J'| \le \theta$, then there is an order-preserving injective map $g$ from $J'$ into $J$). 
    \end{enumerate}

    To finish the proof it suffices to prove: 

    \begin{enumerate}
        \item[$\boxplus$]  $\EM(J) \in K_\theta$ is universal.
    \end{enumerate}

    Why $\boxplus$ holds?  $\EM(J)$ is a model of power $\theta$ since max$(|J|, \chi) \le \theta$ and $\theta = |J| \le \| \EM(J) \|.$ Let us show that $\EM(J)$ is universal. Suppose that $N \in K_\theta$.  Applying Fact \ref{X1.2} to a suitably large $\Op(N)$ find $M \in K_\lambda, N \leq_{\gk} M$, so that by $\lambda$-categoricity of $\gk, M \cong \EM(I)$.  There is a surjective $\leq_{\gk}$-elementary embedding $h\colon  N \underset{\gk}  \rightarrow N' \leq_{\gk} \EM(I)$ and there exists $J' \subseteq I, \,  | J'| \le \| N' \| + \chi = \theta$, such that $N' \subseteq \EM(J')$.  So by $(\ast)$ there is an order preserving injective map $g$ from $J'$ into $J$. Now $g$ induces an $\leq_{\gk}$-elementary
    embedding $\hat g$ from $\EM(J')$ into $\EM(J)$.  Let $f = \hat g h$, then $f\colon  N \underset{\gk} \rightarrow \EM(J)$ is as required.
\end{PROOF}

\begin{Theorem}\label{X3.12}
    Suppose that $\gk$ is categorical in $\lam$, $\theta \in [\k+ |T|, \lam)$, $N\in K_{<\lam}$ is nice, $M\in K_{\theta}$ and $M{\underset{\text{nice}}  \leq} N$. Then $M$ is nice.
\end{Theorem}

\begin{PROOF}{\ref{X3.12}}
    Let $B\in K_{\theta}$, $M \leq_{\gk} B$ and we show that $M \underset{\rm{nice}}{\leq} B$. Well, since $M
   \underset{\rm{nice}}{\leq}
    N$ and $M \leq_{\gk} B$, by \ref{X2.1} there exists an amalgam $M^* \in K_{<\lam}$ of $N, B$ over $M$. Without loss of generality by \ref{X1.5} $\Vert M^*||=||N||$. But $N$ is nice, hence $N \underset{\rm{nice}}{\leq} M^*$. Since $M\prenice N$, it follows by \ref{w17} that $M \underset{\rm{nice}}{\leq} M^*$. Since $M \leq_{\gk} B \leq_{\gk} M^*$, it follows by \ref{X1.11} that
    $M \underset{\rm{nice}}{\leq} B$.
\end{PROOF}

\newpage 

\section{$(\theta, \partial)$-saturated models}\label{4}

In this section, we define  notions of saturation which will be of use in proving
amalgamation for $\gk_\lambda$.

\begin{Definition}\label{X4.1}
    Suppose that $\partial$ is an ordinal,
    $\aleph_{0} \le \partial \le \theta \in [\chi, \lambda)$.
    
    (1) An $\tau$-structure $M$ is $(\theta, \partial)$-saturated\footnote{Called $(\theta, \partial)$-trimmed in \cite{Sh:600}.} iff:
    
    \begin{enumerate}
        \item[(a)] $\| M \| = \theta,$
    
        \item[(b)] there exists a continuous $\leq_{\gk}$-chain $\langle M_i \in K_\theta \colon  i < \partial \rangle$ witnessing it, which means:
    
        \begin{enumerate}
            \item[(i)] $M_0$ is nice and universal, 
        
            \item[(ii)] $M_{i +1}$ is universal over $M_i$, 
        
            \item[(iii)] $M_i$ is nice, and,
            
            \item[(iv)] $M = \underset{i < \partial} \cup M_i.$
        \end{enumerate}
    \end{enumerate}
    
    (2) $M$ is $\theta$-saturated iff $M$ is $(\theta$,  $\cf(\theta))$-saturated.
     
    (3) $M$ is $(\theta, \partial)$-saturated over $N$ iff $M$ is $(\theta, \partial)$-saturated as witnessed by a chain
    $\lng M_i\colon i<\partial\rng$ such that  $N \leq_{\gk} M_0$.
    
    The principal facts established in this section connect the existence, uniqueness and niceness of $(\theta, \partial)$-saturated models.
\end{Definition}

\begin{Theorem}\label{X4.2}
    Suppose that $\gk$ is categorical in $\lambda$ and $\partial \le \theta \in [\chi , \lambda)$.  \underline{Then}:
    
    (1) there exists a $(\theta, \partial)$-saturated model $M,$
    
    (2) for $\partial$ a limit ordinal, $M$ is unique up to isomorphism,
    
    (3) $M$ is nice.
\end{Theorem}

\begin{PROOF}{\ref{X4.2}} 
    One proves (1), (2), and (3)
    simultaneously by induction on $\partial$.
    
    \un{Ad (1)}.  Choose a continuous $\leq_{\gk}$-chain $\langle M_i \in K_\theta\colon  i < \partial \rangle$ of nice models by induction on $i$ as follows.  For $i = 0$, apply \ref{X3.11} to find a nice universal model $M_0 \in K_\theta$.  For $i = j +1$, note that $M_j$ is an a.b. by \ref{X3.5} (since $M_j$ is nice), hence by \ref{X3.10} there exists a nice model $M_i \in K_\theta, M_j \leq_{\gk} M_i, M_i$ universal over $M_j$.  For limit $i$, let $M_i = \underset{j < i} \cup M_j$.  Note that by the inductive hypothesis (3) on $\partial$ for $i < \partial$, since $M_i$ is $(\theta, i)$-saturated, $M_i$ is nice. Thus $M = \underset{i < \partial} \cup M_i$ is $(\theta, \partial)$-saturated (witnessed by $\langle M_i \colon  i < \partial \rangle$).  Note that $M$ is universal since $\langle M_i \colon  i < \partial \rangle$ is continuous and $M_0$ is universal.
    
    \un{Ad (2)}. Recall that each $M_{i}$ is an amalgamation base by $\ref{X2.1}.$  As $\partial$ is a limit ordinal
    standard back-and-forth argument shows that if $M$ and $N$ are $(\theta, \partial)$-saturated models, then $M$ and $N$ are isomorphic.
    
    \un{Ad (3)}.  By the uniqueness (i.e. by Ad(2)) it suffices to prove that some $(\theta, \partial)$-saturated model is nice. Suppose that $M$ is $(\theta, \partial)$-saturated. We'll show that $M$ is nice.
    
    If $ \cf(\partial) < \partial$, then $M$ is also ($\theta$, $\cf(\partial$))-saturated and hence by the inductive hypothesis (3) on $\partial$ for cf$(\partial), M$ is nice. So we'll assume that cf$(\partial) = \partial$.  Choose a continuous $\leq_{\gk}$-chain $\langle M_i \in K_\theta \colon  i < \theta^+ \rangle$ such that: $M_0$ is nice and universal (possible by \ref{X3.11}); if $M_i$ is nice, then $M_{i +1} \in K_\theta$ is nice and universal over $M_i$ (possible by \ref{X3.5} and \ref{X3.10}); if $M_i$ is not nice (so necessarily $i$ is a limit ordinal), then $M_{i +1} \in K_\theta, M_i \leq_{\gk} M_{i +1}$ and $M_i \underset{\text{nice}}  \nleq M_{i +1}$. By \ref{X3.1} and \ref{X1.10} there is a club $C$ of $\theta^+$ such that if $i \in C$, then $M_i \underset{\text{nice}}  \leq M_{i +1}$.  So by the choice of $\langle M_i \colon  i < \theta^+ \rangle$, if $i \in C$, then $M_i$ is nice.  Choose $i \in C, i = \sup (i \cap C)$, $\cf(i) = \partial$. It suffices to show that $M_i$ is $(\theta, \partial)$-saturated (for then by (2) $M_i$ is isomorphic to $M$ and so $M$ is nice). Choose a continuous increasing sequence $\langle \alpha_\zeta\colon  \zeta < \partial \rangle \subseteq C$ such that $i = \bigcup \{ \alpha_{\zeta}\colon  \zeta < \partial \}$ (recall that $i = \sup (i \cap C)$,  $\cf(i) = \partial$).  Now $M_i = \underset{\zeta < \partial} \cup M_{\alpha_\zeta}$.  Of course $M_{\alpha_0}$ is universal (since $M_0$ is universal and $M_0 \leq_{\gk} M_{\alpha_0}$), $M_{\alpha_{\zeta+1}}$ is universal over $M_{\alpha_\zeta}$ since $M_{\alpha_\zeta +1}$ is universal over $M_{\alpha_\zeta}$ and $M_{\alpha_\zeta} \leq_{\gk} M_{\alpha_\zeta +1} \leq_{\gk} M_{\alpha_{\zeta + 1}}$. Also $M_{\alpha_\zeta}$ is nice for each $\zeta < \partial$ since $\alpha_\zeta \in C$.  Hence $M_i$ is $(\theta, \partial)$-saturated, recall that $M_{i}$ is nice because $i \in C,$ so we are done. 
\end{PROOF}

\begin{Remark}\label{X4.3}
    Remember that by \ref{X3.12}, if $\gk$ is categorical in $\lam$, $\theta \in [\chi, \lam)$, $N\in K_{<\lam}$ is nice, $M\in K_{\theta}$ and $M{\underset{\text{nice}}  \leq} N,$ then $M$ is nice.
\end{Remark}

\begin{Theorem}\label{X4.4}
    Suppose that $\gk$ is categorical in $\lam$, $\chi \le \theta<\theta^+<\lam$. If $\lng M_i\in K_{\theta}\colon i< \theta^+\rng$ is a continuous $\leq_{\gk}$-chain of nice models such that $M_{i+1}$ is universal over $M_i$ for $i<\theta^+$, \underline{then} $\bigcup\limits_{i<\theta^+}M_i$ is $(\theta^+, \theta^+)$-saturated.
\end{Theorem}

\begin{remark}\label{X4.4A}
    Why this is not trivial? Because here $M_{i}$ is of cardinality $\theta$ whereas in Definition \ref{X4.1} the $M_{i}$ are of cardinality $\theta^{+}$. 
\end{remark}

\begin{PROOF}{\ref{X4.4}}
    Write $M=\bigcup\limits_{i<\theta^+}M_i$. Note that if $\lng M'_i\in K_{\theta} \colon i<\theta^+ \rng$ is any other continuous $\leq_{\gk}$-chain of nice models such that $M'_{i+1}$ is universal over $M'_i$ then $\bigcup\limits_{i<\theta^+}M'_i \cong M$ (use again the back and forth argument recalling that $M_{0},$ is an a.b., so as $M_{j}$ is universal over $M_{0},$ it is universal).
    
    By Theorem \ref{X4.2} there exists a $(\theta^+,\theta^+)$-saturated model $N$ which is unique and nice. In particular $||N||=\theta^+$ and there exists a continuous $\leq_{\gk}$-chain
    $\lng N_i\in K_{\theta^+}\colon i< \theta^+\rng$ such that:
    
    \begin{enumerate}
        \item[(i)] $N_0$ is nice and universal,
        
        \item[(ii)] $N_{i+1}$ is universal over $N_i$,
        
        \item[(iii)] $N_i$ is nice,
        
        \item[(iv)] $N=\bigcup\limits_{i<\theta^+}N_i$. 
    \end{enumerate}    

    It suffices to prove that $M$ and $N$
    are isomorphic models.
    
    Without loss of generality $|N|=\theta^+$. By Fact \ref{X1.2}, the set $C_1=\{ \d < \theta^+ \colon  N\uhr\delta \leq_{\gk} N \}$ contains a club of $\th^+$. By \ref{X3.1} there exists a club $C_2 \subseteq C_{1}$ of $\th^+$ such that for every $\d\in C_2$, $N\uhr \d \underset{\rm{nice}}{\leq} N$. Since $\{ |N_i|  \colon i < \th^+ \}$ is a continuous increasing sequence of subsets of $\th^+$, it follows that $C_3=\{ \d<\th^+  \colon \d \subseteq |N_{\d}|\}$ is a club of $\th^+$. Hence there is a club $C_4$ of $\th^+$ such that $C_4\subseteq C_1\cap C_2 \cap C_3\cap [\th, \th^+)$. Note that for $\d\in C_4$ one has $N\uhr \d \prenice N$, $|N\uhr \d|= \d \subseteq |N_{\d}|$ and $N_{\d}  \leq_{\gk} N$, so that $N\uhr \d \leq_{\gk} N_{\d} \leq_{\gk} N$ and so by \ref{X1.10} $N\uhr \d \prenice N_\d$. $\lng N_{\d} \colon \d\in C_4\rng$ is a continuous increasing $\prec_{\gk}$-chain, $N_{\d}\in K_{\th^+}$ and $N\uhr \d\in K_{\th}$.
    
    By \ref{X3.12} $N\uhr\d$ is nice since $N_{\d}$ is nice (by (iii)). So by \ref{X3.10} 
    $N\uhr{\d}$ has a nice $\leq_{\gk}$-extension $B_{\d}\in K_{\th}$ which is universal over $N\uhr \d,$ without loss of generality $N\uhr \d \leq_{\gk} B_{\d} \leq_{\gk} N$.
    
    [Why? since $N\uhr \d \leq_{\gk} B_{\d}$ (in fact $N\uhr \d \prenice B_{\d}$) and
    $N\uhr \d\prenice N_{\d}$, by \ref{X2.1} there exists an amalgam $A_{\d}\in K_{\le \th^+}$ of
    $B_{\d}$, $N_{\d}$ over $N\uhr \d$. Let $f_{\d} : B_{\d} \underset{\gk} \to A_{\d}$ be a witness. But $N_{\d+1}$ is universal over $N_{\d}$ (by (ii)), so $A_{\d}$ can be $\leq_{\gk}$-elementarily embedded into $N_{\d+1}$ over $N_{\d}$ (say by $g_{\d}$), hence $B_{\d}$ can be $\leq_{\gk}$-elementarily embedded into $N$ (using $g_{\d}f_{\d}$).]
    
    Let $C_5=\{\d\in C_4\colon \text{if } \al \in C_4\cap \d$, then $|B_{\al}|\subseteq \d\}$.
    Note that $C_5$ is a club of $\th^+$ since $||B_{\al}||=\th$. [For $\al \in C_4$, let $E_{\al}=(\sup |B_{\al}|, \th^+)\cap C_4$, let $E_{\al}=\th^+$ for $\al\not\in C_4$ and let $E$ be the diagonal intersection of $\lng E_{\al} \colon  \al <\th^+\rng$, i.e. $E=\{ \d<\th^+\colon (\forall \al <\d)(\d\in E_{\al})\}$. Note that $E$ is a club of
    $\th^+$ and $C_5 \supseteq E\cap C_4$ which is a club of $\th^+$.] 
    
    Thus $\lng N\uhr\d\colon \d\in C_5\rng$ is a continuous $\leq_{\gk}$-chain of nice models, each of power $\th$. If $\d_1\in C_5$ and $\d_2=\min(C_5\setminus (\d_1+1))$, then
    $N\uhr \d_1 \leq_{\gk} B_{\d_1} \leq_{\gk} N\uhr \d_2$. Hence $N\uhr \d_2$ is universal over $N\uhr\d_1$ (since $B_{\d_1}$ is universal over
    $N\uhr\d_1$). Let $\{\d_i \colon  i<\th^+\}$ enumerate $C_5$ and set $M'_i=N\uhr\d_i$. Note that $N=\bigcup\limits_{i<\th^+}M'_i$. Then $\lng M'_i\in K_{\th}\colon  i< \th^+\rng$
    is a continuous $\leq_{\gk}$-chain of nice models, $M'_{i+1}$ is universal over $M'_i$.
    Therefore $N$ and $M$ are isomorphic
    (as said at the beginning of the proof), as required.
\end{PROOF}

\begin{notation}\label{X4.5}
    $\Theta=\{\ov\th \colon \ov\th =\lng \th_i\colon  i< \d\rng$ is a (strictly) continuous increasing sequence of cardinals, $\chi < \th_0$, $\d<\th_0$ (a limit ordinal), $\bigcup\limits_{i\le\d}\th_i\le\lam\}$ and $\Theta^-=\{\bar\theta\in \Theta\colon \sup\theta_i<\lambda\}$.
\end{notation}

\begin{Remark}\label{X4.6}
    Let $\th = \sup (\ov\th) = \sup \{ \theta_{i} \colon i < \lg(\bar{\theta}) \}$ for $\bar\th\in\Theta$. Then $\th$ is singular, since $\cf(\th)\le \d<\th_0\le\th$.
\end{Remark}

\begin{Definition}\label{X4.7}
    Let $\bar \th \in \Theta$. A model $M$ is $\ov\th$-saturated iff there is a continuous $\leq_{\gk}$-chain $\lng M_i\in K_{\th_i} \colon  i<\d\rng$ such that $M=\bigcup\limits_{i<\d}M_i$, $M_i$ is nice and $M_{i+1}$ is $\th_{i+1}$-universal over $M_i$.
\end{Definition}

\begin{Definition}\label{X4.8}
    Suppose that $\ov\th\in \Theta$. $\Pr(\ov\th)$ holds iff every $\ov\th$-saturated model is nice.
\end{Definition}

\begin{Remark}\label{X4.9}
    (1) If $\bar\theta_1$, $\bar\theta_2\in\Theta$, $\mathrm{rng}(\bar\theta_1)\subseteq \mathrm{rng}(\bar\theta_2)$, $\sup\mathrm{rng}(\bar\theta_1)=\sup\ \mathrm{rng}(\bar\theta_2)$, and $M$ is $\bar\theta_2$-saturated, then $M$ is $\bar\theta_1$-saturated.

    (2) For $\bar\theta\in \Theta^-$ and $\Pr(\bar \theta^\prime)$ whenever $\bar \theta^\prime\in \Theta$ is a proper initial segment of $\bar \theta$, there is a $\bar\theta$-saturated model and it is unique.
\end{Remark}

\begin{Theorem}\label{X4.10}
    Suppose that $\ov\th \in \Theta^-$ and for every limit ordinal $\al<\rm{lg}(\ov\th)$, $\Pr(\ov\th\uhr\al)$. Then $\Pr(\ov\th)$.
\end{Theorem}

\begin{PROOF}{\ref{X4.10}}
    Let $\theta = \sup(\bar{\theta}).$ By \ref{X4.9}(1) and the uniqueness of $\ov\th$-saturated models \ref{X4.9}(2), without loss of generality one may assume that $\lg(\ov\th) = \cf(\sup (\ov\th)) = \cf(\theta)$. Now, by \ref{X4.6}, we know $(\cf(\th))^+<\th \ (= \sup(\bar{\theta}))$, so by \cite[\ref{X1.5} + Fact \ref{X1.2}(1)]{Sh:420} 
    there exists $\lng S, \lng C_{\al} 
 \colon \al \in S\rng\rng$ such that: 
    
    \begin{itemize}
        \itm{$(\al)$} $S\subseteq \th^+$ is a set of ordinals; $ 0 \notin S,$
        
        \itm{$(\b)$} $S_1=\{ \al\in S\colon  \cf(\al)= \cf (\th)\}$ is a stationary subset of $\th^+,$
        
        \itm{$(\gamma)$} if $\alpha \in S$ 
        then $\al=\sup (C_{\al})$ and, if $\alpha\in S \setminus S_{1}$ then $\rm{otp}(C_{\al}) < \cf (\th),$
    
        \itm{$(\d)$} if $\b \in C_{\al}$, then $\beta \in S$ and $C_{\b}=C_{\al}\cap \b,$
    
        \itm{$(\epsilon)$} $C_{\al}$ is a set of successor ordinals.
    \end{itemize}
    
    [Note that the existence of $\lng S,\lng C_{\al} \colon \al \in S\rng\rng$
    is provable in ZFC.] 
    
    Without loss of generality $S\setminus S_1=\cup\{C_\alpha \colon \alpha\in S_1\}$.
    We shall construct the required model by induction, using $\lng  C_{\al} \colon \al \in S\rng$.  Remember $\bar \theta=\langle \theta_\zeta\colon  \zeta< \cf(\theta)\rangle$.
    Let us start by defining by induction on $\al<\th^+$ the following entities: $M_{\al}$, $M_{\al, \xi}$ (for $\alpha<\theta^+$, $\xi < \cf(\th))$, and $N_\al$ (only when $\al\in\bigcup\limits_{\beta\in S}C_\beta$)  such that:
    
    \item{(A1)} $M_{\al} \in K_{\th}$,
    
    \item{(A2)} $\lng M_{\al} \colon  \al < \th^+\rng$ is a continuous increasing $\leq_{\gk}$-chain of models,
    
    \item{(A3)} $M_{\al +1}$ is nice, and if $M_{\al}$ is not nice, then $M_{\al} {\underset  {\text{nice}}  {\not\nleq}}  M_{\al+1},$
    
    \item{(A4)} $M_{\al} \not= M_{\al+1},$
    
    \item{(A5)} $M_{\al+1}$ is weakly universal over $M_{\al},$
    
    \item{(B1)} $M_{\al}=\bigcup\limits_{\xi<\cf(\th)}M_{\al, \xi}$, $||M_{\al, \xi}||=\th_{\xi},$

    \item{(B2)} if $\al \in S_1 \; , \; \b \in C_{\al} \; , \; \gamma\in C_{\al} \; , \; \b < \gamma$, then:

    \begin{enumerate}
        \item[(a)]   $N_{\b} \leq_{\gk} M_{\b}$,

        \item[(b)] $||N_{\b}||=\th_{\otp(C_{\b})}$,

        \item[(c)] $(\forall \xi < \otp(C_{\b}))(M_{\b, \xi} \leq_{\gk} N_{\gamma}),$

        \item[(d)] $N_{\b}$ is nice,

        \item[(e)] $N_{\gamma}$ is $\th_{\otp(C_{\gamma})}$-universal over $N_{\b}$.
    \end{enumerate}
    
    There are now two tasks at hand.
    First of all, we shall explain how to construct these entities (THE CONSTRUCTION, below).
    Then we shall use them to build a nice $\ov\th$-saturated  model (PROVING  $\Pr(\ov\th)$, below). \relax From the uniqueness of $\ov\th$-saturated models it will
    thus follow that $\Pr(\ov\th)$ holds.

    \underline{THE CONSTRUCTION}: we consider several cases: 
    
    \underline{Case (i)}: $\b=0$. Choose $M_0\in K_{\th}$ and $\lng M_{0, \xi} \in K_{\th} \colon  \xi< \cf(\th) \rng$
    with $M_0=\bigcup\limits_{\xi<\cf(\th)}M_{0, \xi}$ using Fact \ref{X1.2}. There is no need to define $N_0$ since $0\not\in C_{\al}$.
    
    Case (ii): $\b$ is a limit ordinal. Let $M_{\b}=\bigcup\limits_{\gamma<\b}M_{\gamma}$
    and choose $\lng M_{\b, \xi} \colon  \xi < \cf(\th)\rng$ using Fact \ref{X1.2}. Again there's no call
    to define $N_{\b}$ since $C_{\al}$ is always a set of successor ordinals.
    
    Case (iii): $\b$ is a successor ordinal, $\b=\gamma+1$. Choose $M'_\gamma\in K_\theta$ such that $M_\gamma \leq_{\gk} M'_{\gamma}$ and if possible $M_\gamma \underset{\text{nice}}  {\nleq} M'_\gamma$; without loss of generality $M'_\gamma$ is weakly universal over $M_\gamma$. If $\b\not\in S$, then define things as above, taking into account (A3).  The definitions of $M_{\b}$, $M_{\b, \xi}$ present no special difficulties. Now suppose that $\b\in S$. The problematic entity to define is $N_{\b}$.

    If $C_{\b}=\emptyset$, choose for $N_{\b}$ any nice sub-model (of power $\theta_{0}$) of
    $M_{\gamma}$.

    If $C_{\b}\not=\emptyset$, then first define $N^-_{\b}=\bigcup\limits_{\gamma\in C_{\b}}N_{\gamma}$. Note that $N^-_{\b}$ is nice. [If $C_{\b}$ has a last element $\b'$, then $N^-_{\b}=N_{\b'}$ which is nice; if $C_{\b}$ has no last element, then $N^-_{\b}=\bigcup\limits_{\gamma\in C_{\b}}N_{\gamma}$ is $\ov\th\uhr \otp(C_{\b})$-saturated, and, by the hypothesis of the theorem, $\Pr(\ov\th\uhr \otp(C_{\b}))$, so $N^-_{\b}$ is nice.] Also $N^-_\beta \preceq_{\gk} M_\gamma$. If $\rm{otp}(C_\b)$ is a limit ordinal we let $N_\b=N^-_\b$ and $M_\beta= M^\prime_\gamma$, so we have finished, so assume $\otp(C_\b)$ is a successor ordinal. To complete the definition of $N_{\b}$, one requires a Lemma (the proof of which is similar to \ref{X3.9}, \ref{X3.10}):
    
    \begin{enumerate}
        \item[{$(\ast)$}] if $A\subseteq M\in K_{\th}$, $|A|\le \th_j<\th$, then there exist a nice $M^+\in K_{\th}$, $M \leq_{\gk} M^+$, and nice models $N^*$, $N^+\in K_{\th_j}$, $A \subseteq N^* \leq_{\gk} N^+ \leq_{\gk} M^+$ and $N^+$ is universal over $N^*$.
    \end{enumerate}
    
    Why is this enough? Use the Lemma with $M=M'_{\b}$ and $A=N^-_{\b}\cup \bigcup\limits_{\scriptstyle \xi< \otp(C_{\b}) \atop \scriptstyle \gamma\in C_{\b}}M_{\gamma, \xi}$  to find $N^*$, $N^+$, $M^+$ and choose $N^+$, $M^+$ as $N_\beta$, $M_\beta$ respectively. 
    
    Now, why $(*)$ holds? The proof of $(\ast)$ is easy as $M'_\beta$ is nice.
    
    \underline{PROVING Pr$(\ov\th)$}:

    For $\al \in S_1$, consider $\lng N_{\b} \colon  \b \in C_{\al}\rng$. For $\b,\gamma\in C_{\al}, \b<\gamma$,  one has by (B2)(c) $\bigcup\limits_{\xi < \otp (C_{\b})}M_{\b, \xi} \subseteq N_{\gamma}$. Therefore $M_{\b} \subseteq \bigcup\limits_{\gamma\in C_{\al}}N_{\gamma}.$ $(M_\b=\bigcup\limits_{\xi<\cf(\th)}M_{\b, \xi}=\bigcup\limits_{\xi<\cf(\al)} M_{\b, \xi} $  (since $\al\in S_1)$; for $\xi <\cf(\al)$, choose $\gamma \in C_{\al}$,
    $\xi<\gamma, \b<\gamma$; so $M_{\b, \xi}\subseteq N_{\gamma}$ and $M_{\b}\subseteq \bigcup\limits_{\gamma\in C_{\al}}N_{\gamma})$.

    Thus for every $\b\in C_{\al}$,
    $M_{\b}\subseteq \bigcup\limits_{\gamma\in C_{\al}}N_{\gamma}$ hence $M_\alpha = \bigcup\limits_{\beta\in C_\alpha} M_\beta \subseteq \bigcup\limits_{\gamma\in C_\alpha} N_\gamma$ (remember $\alpha = \sup (C_\alpha)$ as $\alpha \in S_1$). If $\gamma\in C_{\al}$, then $N_{\gamma} \leq_{\gk} M_{\gamma}$ (by (B2)(a)), and so $\bigcup\limits_{\gamma\in C_{\al}}N_{\gamma} \subseteq \bigcup\limits_{\b\in C_{\al}}
    M_{\b}=M_{\al}$ by continuity.
    So $M_\alpha = \bigcup\limits_{\beta\in C_\alpha} N_\beta$ hence $\langle N_\beta\colon  \beta\in C_\alpha\rangle$ exemplifies $M_\alpha$ is
    $\bar \theta$-saturated (remember $\Pr(\bar \theta \uhr \delta)$ for every limit $\delta< \rm{lg}(\bar \theta)$). So $M_{\al}$ is $\ov\th$-saturated for every $\al\in S_1$. In other words $\{\al<\th^+ \colon  M_{\al}$ is $\ov\th$-saturated$\} \supseteq S_1$ and is stationary, so, applying \ref{X3.1}, there exists
    $\al<\th^+$ such that $M_{\al}$ is $\ov\th$-saturated and $M_{\al}\prenice \bigcup\limits_{\b<\th^+}M_{\b}$.  Hence by \ref{X1.10} $M_{\al} \underset{\rm{nice}}{\leq} M_{\al+1}$ and so, since $M_{\al+1}$ is nice (A3), $M_{\al}$ is nice (by \ref{X3.12}).
    
    We conclude that Pr$(\ov\th)$ holds.
\end{PROOF}

To round off this section of the paper, let us make the connection between $\ov\th$-saturation and $(\th, \cf(\th))$-saturation (Notation follows \ref{X4.5}--4.10).

\begin{Theorem}\label{X4.11}
    Let $\bar\theta\in \Theta^-$ and $\theta= \underset{i} \sup (\theta_i)$. Every $\ov\th$-saturated model is $(\th, \cf(\th))$-saturated.
\end{Theorem}

\begin{PROOF}{\ref{X4.11}}
    Let $\lng M_\al\colon \al<\theta^+\rng$ be as in the
    proof of \ref{X4.10}. By \ref{X3.1} there exists a club $C$ of $\th^+$ such that for every $\al\in C$, $M_{\al} \underset{\rm{nice}}{\leq} \bigcup \limits_{\b< \theta^{+}} M_{\b}$ hence by the construction $M_{\al}$ is nice. So if $\al, \b\in C$ and $\al<\b$, then $M_{\b}$ is a universal  extension of $M_{\al}$ and for $\gamma=\sup(\gamma\cap C)$, $\gamma \in C$, one has that  $M_{\gamma}$ is $(\th, \cf(\gamma))$-saturated. Choose $\gamma\in S_1\cap C$ and $\sup(\gamma\cap C)=\gamma$. So $M_{\gamma}$ is $(\th, \cf(\th))$-saturated
    and also $\bar\theta$-saturated (see proof of \ref{X4.10}). Together we finish.
\end{PROOF}

\newpage

\section{THE AMALGAMATION PROPERTY FOR $\gk_{< \lambda}$}\label{5}

Corollaries \ref{X5.4} and \ref{X5.5} are the goal of this section, showing that if $K$ is categorical in $\lambda$ then every element of $\gk_{< \lambda}$ is nice (see \ref{X5.4}) and $\gk_{< \lambda}$ has the amalgamation property (see \ref{X5.5}).

\begin{Lemma}\label{X5.1}
    Suppose  that $\mu$ is singular, $\lng \mu_i \colon  i< \cf(\mu) \rng$ is a continuous strictly increasing sequence of ordinals, $\mu={\underset{i< \cf(\mu)} \sup} \mu_i$, and $\chi \le \mu_0<\mu\le\lam$. Then there exist a linear order $I $ of power $\mu$ and a continuous increasing sequence $\lng I_i\colon i<\cf(\mu)\rng$ of linear orders such that:
    
    \begin{enumerate}
        \item[(a)] $\chi \le |I_i|\le\mu_i$ and $\vert I_{i} \vert < \vert I_{i + 1} \vert$ for each $i$,  
    
        \item[(b)] $\bigcup\limits_{i<\cf(\mu)}I_i=I$,
        
        \item[(c)] every $t\in I_{i+1}\setminus I_i$ defines a Dedekind cut of $I_i$ in
        which (at least) one side of the cut has cofinality $\k$.
    \end{enumerate}
\end{Lemma}

\begin{PROOF}{\ref{X5.1}}
    Let $I=(\{0\}\times \mu)\cup(\{1\}\times\kappa)$, $I_i=(\{0\}\times \mu_i)\cup (\{1\}\times \k)$ ordered by: $$(i,\al)_{<_I} (j,\b) \ \rm{iff} \ i<j \ \rm{or} \  (0 = i = j \ \rm{and} \ \alpha < \beta) \ \rm{or} \  (1=i=j \ \rm{and} \ \al>\b).$$
\end{PROOF}

\begin{Lemma}\label{X5.2}
    Suppose that $\gk$ is categorical in $\lambda > \cf(\lambda), \, \k + \LST_{\gk} < \mu \le \lambda$. If $M \in K_\lam$, \underline{then} there exists a continuous increasing $\prec_{\gk}$-chain $\langle M_i \colon  i < \cf(\lam) \rangle$ of models such that:
    
    \begin{enumerate}
        \item[$($\rm{a})] $ M \leq_{\gk} \underset{i < \cf(\lam)} \bigcup M_i,$
    
        \item[$($\rm{b})] $ \| \underset{i < \cf(\lam)} \cup M_i \| = \lam,$

        \item[$($\rm{c})] $\k + | T| \le \| M_i \| < \| M_{i +1} \| < \lam,$

        \item[$($\rm{d})] for each $i < \cf(\lam) , M_i \underset{\text{nice}}  \leq \left( \underset{j < \cf(\lam)} \cup M_j \right)$. 
    \end{enumerate}
\end{Lemma}

\begin{PROOF}{\ref{X5.2}}
    As $\lambda$ is a limit cardinal, choose a continuous increasing sequence $\langle \mu_i \colon i < \cf(\lam) \rangle$, $\lam = \underset{i < \cf(\lam)} \sup \mu_i$, $\k +|T| \le \mu_0 < \lambda$. Let $\langle I, \langle I_i \colon i < \cf(\lam) \rangle \rangle$ be as in \ref{X5.1}. By $\lambda$-categoricity of $\gk$  without loss of generality $M = \EM(\lambda)$.  Let $M_i = \EM(I_i)$ for $i < \cf(\lam)$.  Clearly (a), (b), and (c) hold.  To obtain (d), observe that by \ref{X2.6} and \ref{X3.5} it suffices to remark that by demand (c) from \ref{X5.1} on $\langle I_i \colon i < \cf(\lam) \rangle$ clauses $(\aleph)$ or $(\beth)$ in \ref{X2.6}  holds for each $t \in I \setminus I_i$. 
\end{PROOF}

\begin{Theorem}\label{X5.3} 
    For every $\mu \in [ \chi, \lambda]$ and $M \in K_\mu$, there exists $M' \in K_\mu, 
    M \leq_{\gk} M'$ such that: 
    
    \begin{enumerate}
        \item[$(\ast)_{M'}$]  for every $A \subseteq | M'|$, $|A|<\lambda\ \wedge\ |A|\leq \mu$, there is $N \in K_{ \chi +|A|}$ such that $A \subseteq N \leq_{\gk} M'$ and $N$ is nice.
    \end{enumerate}
\end{Theorem}

\begin{PROOF}{\ref{X5.3}}  
    The proof is by induction on $\mu$.

    \underline{Case 1}:  $\mu = \chi$.  By \ref{X3.3} there is $M' \in K_\mu, \,  
    M \leq_{\gk} M'$ and 
    $M'$ is nice.  Given $A \subseteq | M'|$ let $N = M'$ and note that $N$ is as required in $(\ast)_{M'}$. 
    
    \underline{Case 2}:   $\chi < \mu$. Without loss of generality, one can replace $M$ by any $\leq_{\gk}$-extension in $K_\mu$. Choose a continuous increasing sequence $\langle \mu_i \colon i < \cf(\mu) \rangle$ such that if $\mu$ is a limit cardinal it is a strictly increasing
    sequence with limit $\mu$; if $\mu$ is a successor, use $\mu_i^+=\mu$ for every $i < \cf(\mu),$ and in both cases
    $\chi \le \mu_i < \mu$.  Find
    $\bar M=\langle M_i \colon i < \cf(\mu) \rangle$ such that:
    
    \begin{enumerate}
        \itm{(a)} $M \leq_{\gk} \bigcup\limits_{i< \cf(\mu)}M_i,$
        
        \itm{(b)} $||\bigcup\limits_{i< \cf(\mu)}M_i||=\mu,$
        
        \itm{(c)} $||M_i||=\mu_i,$
        
        \itm{(d)} $M_i
        \underset{\rm{nice}}{\leq}
        \bigcup\limits_{j< \cf(\mu)}M_j$.
    \end{enumerate}
    
    Why does $\bar M$ exist? If $\mu=\lam$ by \ref{X5.2}, otherwise by \ref{X4.4} ($\mu$ regular) and \ref{X4.11} ($\mu$ singular).

    Choose by induction on $i < \cf(\mu)$ models $L^0_i, L^1_i, L^2_i$ in that order such that:
    
    \begin{itemize}
        \itm{$(\alpha)$} $M_i \leq_{\gk} L^0_i \leq_{\gk} L^1_i \leq_{\gk} L^2_i \in K_{\mu_i},$ 
        
        \itm{$(\beta)$} $j<i \Rightarrow L^2_j \leq_{\gk} L_{i}^0,$

        \itm{$(\gamma)$} $(\ast)_{L^1_i}$ holds , i.e. for each $A \subseteq |L^1_i|$,  there is $N \in K_{\le \k + |T| + |A|}$ such that
        $A \subseteq N \leq_{\gk} L^1_i$ and $N$ is nice (so in particular $L^1_i$ is nice, letting $A = |L^1_i|$),
        
        \itm{$(\delta)$} $L^2_i$ is nice and $\mu_i$-universal over $L^1_i$,
        
        \itm{$(\varp)$} $L^0_i$ is $\leq_{\gk}$-increasing continuous,

        \item[$(\zeta)$] $L_{i}^{\ell} \cap \bigcup_{j < \cf(\mu)} M_{j} = M_{i}$ (or use system of $\leq_{\gk}$-embeddings). 
\end{itemize}

    For $i = 0,$ let $L_{i}^{0} = M_{0}.$ For $i = j +1,$ note that by \ref{X2.1} there is an amalgam $L_{i}^{0} \in K_{\mu_{i}}$ of $M_{i},$ $L_{j}^{2}$ over $M_{j}$ since $M_{j} \underset{\mathrm{nice}}{\leq} M_{i}$ and $M_{j} \leq_{\gk} L_{j}^{2}$ (use last phrase of Fact \ref{X1.2} for clause $(\zeta)$); actually not really needed. For limit $i,$ continuity necessitates choosing $L_{i}^{0} = \bigcup_{j < i} L_{j}^{0}$ (note that in this case $L_{i}^{0} = \bigcup_{j < i}L_{j}^{2}$). To choose $L_{i}^{1}$ apply the inductive hypothesis with respect to $\mu_{i}, L_{i}^{0}$ to find $L_{i}^{1}$ so that $L_{i}^{0} \leq_{\gk} L_{i}^{1}$ and $(\gamma)(\ast)_{(L_{i}^{1})}$ holds. To choose $L_{i}^{2}$ apply \ref{X3.7} to $L_{i}^{1} \in K_{\mu_{i}}$ giving $L_{i}^{1} \leq_{\gk} L_{i}^{2},$ $L_{i}^{2}$ is nice and $\mu_{i}$-universal over $L_{i}^{1}$ (so $(\delta)$ holds).  

    Let $L = \bigcup_{i < \cf(\mu)} L_{i}^{0} = \bigcup_{i < \cf(\mu)} L_{i}^{1} = \bigcup_{i < \cf(\mu)} L_{i}^{2},$ and let $   
    L  
    _i=   
    L
    _i^0$ if $i$ is a  limit,  $ 
    L 
    _i^1$
    otherwise. Now show by induction on $i < \cf(\mu)$ that 
    $  
    L
    _i$ is nice.
    
    [Why? show by induction on $i$ for $i=0$ or $i$ successor 
    that $  
    L
    _i =   
    L
    ^1_i$
    hence use clause $(\gamma)$, if $i$ is limit 
    then $  
    L
    _i$ is $(\bar \theta \uhr i)$-saturated, 
    hence $  
    L
    _i$ is nice by \ref{X4.8}, \ref{X4.10}.]
    
    Now $\langle  
    L
    _i \colon i < \cf(\mu) \rangle$ witnesses that \underbar{if} $\mu$ is regular, $L$ is $(\mu,\mu)$-saturated by \ref{X4.4}, \underbar{if} $\mu$ is singular, $L$ is $\bar\mu$-saturated; in all cases $L$ is $\bar\mu$-saturated of power $\mu$, hence by the results of section 4 (i.e. \ref{X4.8}, \ref{X4.10}) if $\mu<\lam$ then $L$ is nice. Claim \ref{X5.3.1} below provides the desired model $M',$ so we are done. 
\end{PROOF}

\begin{claim}\label{X5.3.1} 
    $M' = L$ is as required.
\end{claim}

\begin{PROOF}{\ref{X5.3.1}}
    $M \leq_{\gk} \underset{i < \cf(\mu)}  \cup M_i  \leq_{\gk} \underset{i < \cf(\mu)}  \cup L^0_i = L \in K_\mu$. Suppose that $A \subseteq |L|$.  If $|A| = \mu$, then necessarily $\mu<\lam$ and we take $N = L$.  So without loss of generality, $|A| < \mu$.  If $\mu = \cf(\mu)$ or $|A| < \cf(\mu)$, then there is $i < \cf(\mu)$ such that $A \subseteq L^1_i$ and, by $(\gamma), (\ast)_{L^1_i}$ holds, so there is $N \in K_{\k + \LST(\gk) + |A|}, A \subseteq N \leq_{\gk} L^1_i, N$ is nice and $N \leq_{\gk} L$ as required.  So suppose that $\cf(\mu) \le |A| < \mu$.  Choose by induction on $i < \cf(\mu)$ models $N^0_i, N^1_i, N^2_i$ in that order such that:
    
    \begin{enumerate}
        \item[$(\alpha)$] $N_{i}^{0} \leq_{\gk} N_{i}^{1} \leq_{\gk} N_{i}^{2},$
        
        \item[$(\beta)$] $N_{i}^{2} \leq_{\gk} N_{i+1}^{0},$
        
        \item[$(\gamma)$] $A \cap L_{i}^{0} \subseteq N_{i}^{0} \leq_{\gk} L_{i}^{0},$  
        
        \item[$(\delta)$] $N^1_i \leq_{\gk} L^1_i$ and $ N^1_i$ is nice,
        
        \item[$(\varp)$] $N^2_i \leq_{\gk} L^2_i, N^2_i$  is nice and  universal over $N^1_i,$
        
        \item[$(\zeta)$] $N^0_i, N^1_i, N^2_i$ have power at most $\min\{\chi + |A|, \mu_i\}$.
    \end{enumerate}

    For $i =0$, apply Fact \ref{X1.2} for $A \cap L^0_0, L^0_0$; for $i = j +1$, apply Fact \ref{X1.2} to find $N^0_i \in K_{\mu_i}, (A \cap L^0_i) \cup N^2_j \subset N^0_i \leq_{\gk} L^0_i$ (in particular $N^2_j \leq_{\gk} N^0_i)$; for limit $i, N^0_i = \underset{j < i} \cup N^0_j$.  To choose $N^1_i$, use $(\ast)_{L^1_i}$ for the set $A_i = N^0_i$ to find a nice $N^1_i \in K_{\le \chi + |A|}, N^0_i \leq_{\gk} N^1_i \leq_{\gk} L^1_i$. Note that $\| N^1_i \| \le \mu_i$.  Finally to choose $N^2_i$ note that by Lemma \ref{X3.9} the model $N^1_i$ has a nice extension $N^+_i$ (of power $\| N^1_i \|$) weakly universal over $N^1_i$.  Now $N^1_i$ is nice, hence $N_{i}^{+}$ is universal over $N^1_i$ (by \ref{X3.6}A(5)) and by  Lemma \ref{X2.1} there is an amalgam $N_i$ of $N^+_i, L^1_i$ over $N^1_i$ such that $\| N_i \| \le \mu_i$.  Since $L^2_i$ is universal over $L^1_i$ one can find an $\leq_{\gk}$-elementary sub-model $N^2_i$ of $L^2_i$ isomorphic to $N_i$. Let $N_i$ be $N_i^0$ if $i$ is a limit, $N_i^1$ otherwise; prove by induction on $i$ that $N_i$ is nice (by Theorem \ref{X4.2}).
    
    Now $\underset{i < \cf(\mu)}  \cup N^0_i$ is an $\leq_{\gk}$-elementary sub-model of $L$ of power at most $\k + |T| + |A|$, including $A$ (by $(\gamma)$) and $\underset{i < \cf(\mu)}  \cup N^0_i$ is $(\chi + |A|, \cf(\mu))$-saturated, hence (by 
    Theorem \ref{X4.2}) nice, as required.
\end{PROOF}

\begin{Corollary}\label{X5.4}
    If $K$ is categorical in $\lambda$  then every element of $K_{< \lambda}$ is nice.
\end{Corollary}

\begin{PROOF}{\ref{X5.4}}
    Suppose otherwise and let $N_0 \in K_{< \lambda}$ be a model which is not nice.
    Choose a suitable $\Op$ such that $\| \Op( N_0 ) \| \ge \lambda$ and by Fact \ref{X1.2} find
    $M_0 \in K_\lambda, N_0 \leq_{\gk} M_0 \leq_{\gk} \Op (N_0)$ i.e. 
    $N_0 
    \underset{\rm{nice}}{\leq}
    M_0$.  It follows that:
    
    \begin{enumerate}
        \item[$\boxplus$] if $N_0 \leq_{\gk} N \leq_{\gk} M_0$ and $N \in K_{< \lambda}$ then $N$ is not nice.
    \end{enumerate}
    
    [Why? By \ref{X4.3}; alternatively,  suppose contrariwise that $N$ is  
    nice. So there is $N_1 \in K_{< \lambda}, N_0 \leq_{\gk} N_1, N_0 
    \underset{\rm{nice}}{\leq}
    N_1,$ so $
    N_0 \underset{\text{nice}}  \leq N$ since 
    $N_0 \underset{\text{nice}}  \leq M_0$ and $N \leq_{\gk} M_0$,  
    hence there is an amalgam $N' \in K_{< \lambda}$ of $N_1, N$ over 
    $N_0$. But $N$ is nice, so $N \underset{\text{nice}}  \leq N'; 
    N_0 \underset{\text{nice}}  \leq N$, so $N_0 \underset{\text{nice}} 
    \leq N'$ and so $N_0 
    \underset{\rm{nice}}{\leq}
    N_1$ contradiction.]
    
    On the other hand, applying \ref{X5.3} for $\mu = \lambda$ there  exists $M' \in K_\lambda$ satisfying $(\ast)_{M'}$.  By $\lambda$-categoricity of $\gk$ without loss of generality, $(\ast)_{M_0}$ holds (see \ref{X5.3}) and $A = |N_0|$ yields a nice model $N \in K_{\k + |T| + \| N_0 \|}$ such that $N_0  \leq_{\gk} N \leq_{\gk} M_0$ contradicting $\boxplus$.
\end{PROOF}

\begin{Corollary}\label{X5.5} 
    If $K$ is categorical in $\lambda$,  then $\gk_{< \lambda}$ has the amalgamation property.
\end{Corollary}

\begin{PROOF}{\ref{X5.5}} 
    As every nice $M \in K_{< \lambda}$ is an amalgamation base (by \ref{X3.5}) we are done by the previous corollary.
\end{PROOF}

\bibliographystyle{amsalpha}
\bibliography{shlhetal}
    
\end{document}